\documentclass[12pt,a4paper]{myarticle}
\newif\ifarxiv
\arxivtrue
\ifarxiv\else
\usepackage[width=14truecm,height=23truecm]{crop}
\advance\hoffset-3.5truecm
\advance\voffset-3.8truecm

\fi

\usepackage{graphicx,color}
\graphicspath{ {images/} }
\usepackage{amsmath,amsfonts,amssymb}
\usepackage{mathrsfs}
\usepackage[T1]{fontenc}
\usepackage[latin1]{inputenc}   
\usepackage[all]{xy}
\SelectTips {cm}{}

\def\varxymatrixc#1#2#3{\vcenter{\kern-#1\hbox{$\displaystyle\xymatrix#2{#3}$}\kern#1}}
\def\xylbl[#1]#2{\ar @{} [#1] |{\hbox{#2}}}

\let\mySfragile\S
\def\myS{\protect\mySfragile}

\usepackage{euscript}

\DeclareMathAlphabet\EuScript{U}{eus}{m}{n}
\SetMathAlphabet\EuScript{bold}{U}{eus}{b}{n}
\let\matheu\EuScript
\DeclareFontFamily{OT1}{bfit}{}
\DeclareFontShape{OT1}{bfit}{m}{n}{<->cmbxti10}{}
\DeclareMathAlphabet {\mathbfit} {OT1} {bfit} {m} {n}

\DeclareFontFamily{OT1}{bfsl}{}
\DeclareFontShape{OT1}{bfsl}{m}{n}{<->cmbxsl10}{}
\DeclareMathAlphabet {\mathbfsl} {OT1} {bfsl} {m} {n}

\DeclareFontFamily{OT1}{smallcaps}{}
\DeclareFontShape{OT1}{smallcaps}{m}{n}{<->cmcsc10}{}
\DeclareMathAlphabet {\mathsc} {OT1} {smallcaps} {m} {n}

\input myfleches.sty

\pdfoutput=1

\let\resetcolor\defaultcolor

\def\gobnextchar#1{}

\def\makeRGB#1#2{\expandafter\edef\expandafter\macroname
\expandafter{\expandafter\gobnextchar\string#1}\expandafter\expandafter\expandafter
\let\expandafter\csname\macroname ColorOrg\endcsname#1
\edef#1{{\noexpand\pdfliteral{#2 rg #2 RG}\noexpand\csname
\macroname ColorOrg\noexpand\endcsname\noexpand\resetcolor}}}

\def\makeRGBOne#1#2{\expandafter\edef\expandafter\macroname
\expandafter{\expandafter\gobnextchar\string#1}\expandafter\expandafter\expandafter
\let\expandafter\csname\macroname ColorOrg\endcsname#1
\edef#1##1{{\noexpand\pdfliteral{#2 rg #2 RG}\noexpand\csname
\macroname ColorOrg\noexpand\endcsname{##1}\noexpand\resetcolor}}}

\def\makeRGBTwo#1#2{\expandafter\edef\expandafter\macroname
\expandafter{\expandafter\gobnextchar\string#1}\expandafter\expandafter\expandafter
\let\expandafter\csname\macroname ColorOrg\endcsname#1
\edef#1##1##2{{\noexpand\pdfliteral{#2 rg #2 RG}\noexpand\csname
\macroname ColorOrg\noexpand\endcsname{##1}{##2}\noexpand\resetcolor}}}

\def\defrobuste#1{{\let\protect0\expandafter\let\csname#1fragile\endcsname0\relax
\expandafter\xdef\csname#1\endcsname{{\protect\csname#1fragile\endcsname}}}}
\def\defRobuste#1{{\let\protect0\expandafter\let\csname#1fragile\endcsname0\relax
\expandafter\xdef\csname#1\endcsname{\protect\csname#1fragile\endcsname}}}
\def\multiexpand#1#2{\spoolmacro{#1}#2\end}
\def\myend{\end}
\def\spoolmacro#1#2{\def\donne{#2}\ifx\donne\myend\let\mnext\relax
\else\csname#1\endcsname{#2}\def\mnext{\spoolmacro{#1}}\fi\mnext}
\def\EXTRA#1{\expandafter\def\csname#1#1#1\endcsname{{\mathord{\matheu #1}}}}
\multiexpand{EXTRA}{ABCDEFGHIJKLMNOPQRSTUVWXYZ}
\def\gras#1{\expandafter\def\csname#1g\endcsname{{\mathord{\mathbfsl #1}}}}
\multiexpand{gras}{ABCDEFGHIJKLMNOPQRSTUVWXYZabcdefhijklmnopqrstuvwxyz}
\def\gothique#1{\expandafter\def\csname#1goth\endcsname{{\mathord{\mathfrak #1}}}}
\multiexpand{gothique}{ABCDEFGHIJKLMNOPQRSTUVWXYZabcdefghijklmnopqrstuvwxyz}
\def\defbbb#1{\expandafter\def\csname#1#1\endcsname{{\mathord{\mathbb#1}}}}
\multiexpand{defbbb}{ABCDEFGHIJKLMNOPQRSTUVWXYZ}
\def\defcal#1{\expandafter\def\csname #1\endcsname{{\mathord{\mathcal #1}}}}
\multiexpand{defcal}{ABCDEFGHIJKLMNOPQRSTUVWXYZ}
\def\defrsfs#1{\expandafter\def\csname #1s\endcsname{{\mathord{\mathrsfs #1}}}}
\multiexpand{defrsfs}{ABCDEFGHIJKLMNOPQRSTUVWXYZ}
\def\defwtilde#1{\expandafter\def
\csname w#1fragile\endcsname{\ifmmode{\widetilde{#1}}\else
$\csname w#1fragile\endcsname$\fi}\defrobuste{w#1}\relax}
\multiexpand{defwtilde}{ABCDEFGHIJKLMNOPQRSTUVWXYZ}
\def\defwtildegras#1{\expandafter\def
\csname w#1gfragile\endcsname{\ifmmode{\widetilde{\mathbfsl#1}}\else
$\csname w#1gfragile\endcsname$\fi}\defrobuste{w#1g}\relax}
\multiexpand{defwtildegras}{ABCDEFGHIJKLMNOPQRSTUVWXYZabcdefghijklmnopqrstuvwxyz}
\def\deftildegras#1{\expandafter\def
\csname t#1gfragile\endcsname{\ifmmode{\tilde{\mathbfsl#1}}\else
$\csname t#1gfragile\endcsname$\fi}\defrobuste{t#1g}\relax}
\multiexpand{deftildegras}{ABCDEFGHIJKLMNOPQRSTUVWXYZabcdefghijklmnopqrstuvwxyz}

\def\mou{plus 0.5pt minus 0.2pt}

\def\mmou{plus -0.8pt minus -0.2pt}

\def\defseccntformat{\def\@seccntformat##1{\csname the##1\endcsname.\ }}
\def\defseccntformatline{\def\@seccntformat##1{\csname the##1\endcsname.}}
\def\defseccntformatEmpty{\def\@seccntformat##1{\csname the##1\endcsname}}

\newskip\intertitleskip
\newskip\secskip
\def\testtitle#1{\if@nobreak
\def\nexttesttitle{\vskip-\lastskip\nobreak\vskip\intertitleskip\nobreak}\relax
\else\def\nexttesttitle{#1}\fi\nexttesttitle}

\def\testtitleline#1{\if@nobreak
\def\nexttesttitle{\leavevmode\vskip-\lastskip\nobreak \vskip-1.7\baselineskip\nobreak}\relax
\else\def\nexttesttitle{#1}\fi\nexttesttitle}

\def\@xsect#1{%
  \@tempskipa #1\relax
  \ifdim \@tempskipa>\z@
    \par \nobreak
    \advance\@tempskipa-\parskip
    \vskip \@tempskipa
    \@afterheading
  \else
    \@nobreakfalse
    \global\@noskipsectrue
    \everypar{%
      \if@noskipsec
        \global\@noskipsecfalse
       {\setbox\z@\lastbox}%
        \clubpenalty\@M
        \begingroup \@svsechd \endgroup
        \unskip
        \@tempskipa #1\relax
        \hskip -\@tempskipa
      \else
        \clubpenalty \@clubpenalty
        \everypar{}%
      \fi}%
  \fi
  \ignorespaces}

\renewcommand\section{\penalty-2500\defseccntformat\@startsection {section}{1}{0pt}%
                                   {-2.5ex \mmou}%
                                   {1ex \mou}%
                                   {\raggedright\normalfont\large\bfseries\boldmath\baselineskip18pt}}
\newcommand\sectionline{\penalty-2500\defseccntformatline\@startsection {section}{1}{0pt}%
                                   {-2.5ex \mmou}%
                                   {-1.ex \mmou}%
                                   {\raggedright\normalfont\large\bfseries\boldmath}}
\renewcommand\subsection{\testtitle{\penalty-1500}\defseccntformat\@startsection{subsection}{2}{\z@}%
                                     {-1em\mmou}%
                                     {0.5ex \mou}%
                                     {\normalfont \normalsize\bfseries\boldmath}}
\newcommand\subsectionline{\testtitle{\penalty-1500}\defseccntformat\@startsection{subsection}{2}{\z@}%
                                     {-1em\mmou}%
                                     {-0.5ex \mou}%
                                     {\normalfont \normalsize\bfseries\boldmath}}
\newcommand\subsectionpar{\testtitle{\penalty-1500}\defseccntformat\@startsection{subsection}{2}{\z@}%
                                     {-1em\mmou}%
                                     {0.5ex \mou}%
                                     {\normalfont \normalsize\bfseries\boldmath}}
\renewcommand\subsubsection{\testtitle{\penalty-1500}\defseccntformat\@startsection{subsubsection}{3}{\z@}%
                                     {-\secskip}%
                                     {0.5ex \mou}%
                                     {\normalfont\normalsize\bfseries\boldmath}}
\newcommand\subsubsectionline{\testtitle{\penalty-1500}\defseccntformatline\@startsection{subsubsection}{3}{\z@}%
                                     {-\secskip}%
                                     {-1.ex \mmou}%
                                     {\normalfont\normalsize\bfseries\boldmath}}
\renewcommand\paragraph{\testtitle\relax\defseccntformat\@startsection{paragraph}{4}{\z@}%
                                    {1.ex \mou}%
                                    {-1ex\mmou}%
                                    {\normalfont\normalsize\bfseries\boldmath}}
\renewcommand\subparagraph{\testtitle\relax\defseccntformat\@startsection{subparagraph}{5}{\parindent}%
                                       {1ex \mou}%
                                       {-1ex\mmou}%
                                      {\normalfont\normalsize\bfseries\boldmath}}

\renewcommand\tableofcontents{%
    \section*{\contentsname
        \@mkboth{%
           \MakeUppercase\contentsname}{\MakeUppercase\contentsname}}%
\begingroup
\vskip0.5ex\parskip0pt    
\@starttoc{toc}%
\endgroup
    }

\renewcommand*\l@part[2]{%
  \ifnum \c@tocdepth >-2\relax
    \addpenalty\@secpenalty
    \addvspace{2.25em \@plus\p@}%
    \setlength\@tempdima{3em}%
    \begingroup
      \parindent \z@ \rightskip \@pnumwidth
      \parfillskip -\@pnumwidth
      {\leavevmode
       \large \bfseries #1\hfil \hb@xt@\@pnumwidth{\hss #2}}\par
       \nobreak
       \if@compatibility
         \global\@nobreaktrue
         \everypar{\global\@nobreakfalse\everypar{}}%
      \fi
    \endgroup
  \fi}
\renewcommand*\l@section[2]{%
  \ifnum \c@tocdepth >\z@
    \addpenalty\@secpenalty
    \setlength\@tempdima{1.5em}%
    \begingroup
      \parindent \z@ \rightskip \@pnumwidth
      \parfillskip -\@pnumwidth
      \leavevmode \bfseries\boldmath
      \advance\leftskip\@tempdima
      \hskip -\leftskip
      #1\nobreak\hfil \nobreak\hb@xt@\@pnumwidth{\hss #2}\par
    \endgroup
  \fi}
\renewcommand*\l@subsection{\@dottedtocline{2}{1.5em}{2em}}
\renewcommand*\l@subsubsection{\@dottedtocline{3}{3.5em}{2.7em}}
\renewcommand*\l@paragraph{\@dottedtocline{4}{7.0em}{4.1em}}
\renewcommand*\l@subparagraph{\@dottedtocline{5}{10em}{5em}}

\newif\ifaddtocontents
\addtocontentstrue
                                      
\def\@sect#1#2#3#4#5#6[#7]#8{%
  \ifnum #2>\c@secnumdepth
    \let\@svsec\@empty
  \else
    \refstepcounter{#1}%
    \protected@edef\@svsec{\@seccntformat{#1}\relax}%
  \fi
  \@tempskipa #5\relax
  \ifdim \@tempskipa>\z@
    \begingroup
      #6{%
        \@hangfrom{\hskip #3\relax\@svsec}%
          \interlinepenalty \@M 
          #8\@@par}%
    \endgroup
    \csname #1mark\endcsname{#7}%
\ifaddtocontents    \addcontentsline{toc}{#1}{%
      \ifnum #2>\c@secnumdepth \else
        \protect\numberline{\csname the#1\endcsname}%
      \fi
      #7}\fi%
  \else
  \ifaddtocontents    
    \def\@svsechd{%
      #6{\hskip #3\relax
      \@svsec #8}%
      \csname #1mark\endcsname{#7}%
        \addcontentsline{toc}{#1}{%
        \ifnum #2>\c@secnumdepth \else
          \protect\numberline{\csname the#1\endcsname}%
        \fi
        #7}}\else
    \def\@svsechd{%
      #6{\hskip #3\relax
      \@svsec  #8}%
      \csname #1mark\endcsname{#7}}\fi%
  \fi
  \@xsect{#5}}

\renewcommand\theenumi{\@arabic\c@enumi}
\renewcommand\theenumii{\@alph\c@enumii}
\renewcommand\theenumiii{\@roman\c@enumiii}
\renewcommand\theenumiv{\@Alph\c@enumiv}

\def\xmygetone#1{\expandafter\mygetone#1}
\def\xmygettwo#1{\expandafter\mygettwo #1}
\def\mygetone#1#2{#1}
\def\mygettwo#1#2{#2}
\def\mygetonetwo#1#2{\def\myrefone{#1}\def\myreftwo{#2}}
\def\mygetrefs#1{\expandafter\let\expandafter\tempref\csname r@#1\endcsname
\expandafter\mygetonetwo\tempref}
\newcounter{myequation}
\def\themyequation{{\rm(\arabic{myequation})}}
\def\eqnnb{\refstepcounter{myequation}\eqno\themyequation}

\def\eqnnbllapraise#1#2{\refstepcounter{myequation}\eqno\hskip0.1pt
\smash{\raise#1\hbox to0pt{\hss\themyequation\label{#2}}}}
\def\eqnollapraise#1#2{\eqno\hskip0.01pt\raise#1\hbox to0pt{\hss$#2$}}

\newif\ifTheoNb
\newif\ifendpoint\endpointtrue
\def\noendpoint{\global\endpointfalse}
\def\endpoint{\ifendpoint.\else\global\endpointtrue\fi}

\def\getTheoName#1#2[#3]{\doTheoHead{#1 #3}{#2}}
\def\doTheoHead#1#2{{\bfseries\ifTheoNb\thesubsubsection. \global\TheoNbfalse\fi#1\endpoint\ }#2\ignorespaces}
\def\makeTheoHead#1#2{\@ifnextchar[{\getTheoName{#1}{#2}}{\doTheoHead{#1}{#2}}}

\newcount\theopenalty
\def\skiptheo{\ifhmode\par\fi
\vskip\secskip\penalty\theopenalty}

\newif\ifenglish

\def\Mynewtheorem#1#2#3{%
\newenvironment{#1}{\intertitleskip-0.2ex\mou\testtitle\skiptheo\relax
\refstepcounter{subsubsection}
\noindent\global\TheoNbtrue\makeTheoHead{#2}{#3}}
{\vskip-\lastskip\nobreak\vskip3pt\mou\global\@nobreakfalse}
\newenvironment{#1*}{\intertitleskip0ex\mou\testtitle\skiptheo\relax
\noindent\global\TheoNbfalse\makeTheoHead{#2}{#3}}
{\vskip-\lastskip\nobreak\vskip3pt\mou\global\@nobreakfalse}}

\def\englishnames{
\def\Conjecturename{Conjecture}
\def\Theoremename{Theorem}
\def\Corollaryname{Corollary}
\def\Propositionname{Proposition}
\def\Lemmaname{Lemma}
\def\Exercisename{Exercise}
\def\Commentname{Comment}
\def\Examplename{Example}
\def\Notationname{Notation}
\def\Definitionname{Definition}
\def\Remarkname{Remark}
\def\vartext{}
\def\defaultdemoname{\gdef\demoname{{\it Proof\endpoint\ }}}
\defaultdemoname}
\englishnames

\Mynewtheorem{conj}{\Conjecturename}{\itshape}
\Mynewtheorem{theo}{\Theoremename}{\itshape}
\Mynewtheorem{coro}{\Corollaryname}{\itshape}
\Mynewtheorem{prop}{\Propositionname}{\itshape}
\Mynewtheorem{lemm}{\Lemmaname}{\itshape}
\Mynewtheorem{exer}{\Exercisename}{\itshape}
\Mynewtheorem{comm}{\Commentname}{}
\Mynewtheorem{comms}{\Commentname s}{}
\Mynewtheorem{exa}{\Examplename}{}
\Mynewtheorem{exas}{\Examplename s}{}
\Mynewtheorem{nota}{\Notationname}{}
\Mynewtheorem{notas}{\Notationname s}{}
\Mynewtheorem{defi}{\Definitionname}{}
\Mynewtheorem{defis}{\Definitionname s}{}
\Mynewtheorem{rema}{\Remarkname}{}
\Mynewtheorem{remas}{\Remarkname s}{}
\Mynewtheorem{var}{\vartext}{\gdef\vartext{}\itshape}

{\catcode`\@=0
\catcode`\{=12\catcode`\}=12
\catcode`\<=1\catcode`\>=2
\catcode`\\=12
@gdef@vartext<@tt@color<red>\def\vartext{...}>>

\def\aliastheorem#1#2{%
\expandafter\def\csname b#1\endcsname##1.{\def\donne{##1}\ifx\donne\rien\begin{#2}\else\begin{#2}[##1]\fi}
\expandafter\def\csname e#1\endcsname{\end{#2}}}

\aliastheorem{p}{prop}
\aliastheorem{t}{theo}
\aliastheorem{l}{lemm}
\aliastheorem{ls}{lemm*}
\aliastheorem{c}{coro}

\newif\ifinsidedemo
\newenvironment{demo}{\ifhmode\par\fi
\global\insidedemotrue
\everypar{}\noindent{\demoname\defaultdemoname}\ignorespaces}{\enddemobox\global\@nobreakfalse\vskip0.5ex \mou}
\def\nblist{\par\nobreak\@nobreaktrue}

\def\enddemobox{\ifinsidedemo
\ifmmode\hbox{$\square$}\else\ifhmode\unskip\else\noindent\fi
\nobreak\null\nobreak\hfill
\nobreak$\square$\fi\fi\global\insidedemofalse}

\def\Jot{0.5\jot}

\let\mathalign\diagalign

\def\mathaligncases#1{\left\{\vcenter{\kern0pt\hbox{$\displaystyle\mathalign{#1}$}\kern0pt}\right.}

\let\dans\subseteq
\let\cont\supseteq
\let\moins\backslash

\def\operator#1#2{\expandafter\def\csname#2fragile\endcsname{\mathop{#1{#2}}\nolimits}\relax
\defRobuste{#2}\relax}
\def\ordinaire#1#2{\expandafter\def\csname#2fragile\endcsname{\mathord{#1{#2}}}\relax
\defRobuste{#2}}
\operator\rm{GL}
\operator\rm{Stab}
\operator\rm{Hom}
\operator\rm{Card}
\operator\rm{dim}
\operator\rm{ker}
\operator\rm{coker}
\operator\rm{Gr}
\operator\rm{id}
\operator\rm{Lie}
\ordinaire\rm{pr}
\operator\rm{im}
\operator\rm{coim}
\operator\bf{Ab}
\operator\rm{Adg}
\operator\rm{Alg}
\operator\rm{Annul}
\operator\rm{Appl}
\operator\rm{Aut}
\operator\rm{Diff}
\operator\rm{End}
\operator\rm{Ext}
\operator\rm{Extgr}
\operator\rm{Torgr}
\operator\rm{Fonct}
\operator\rm{Iso}
\operator\rm{Mod}
\operator\rm{Mor}
\operator\rm{Mor}
\operator\rm{Ob}
\operator\rm{Ouv}
\operator\rm{Parties}
\operator\goth{Rad}
\operator\goth{Rec}
\operator\rm{Red}
\operator\rm{Spec}
\operator\rm{Spm}
\operator\rm{Top}
\operator\rm{Tor}
\operator\rm{Vec}
\def\cl#1{\,\nobreak\overline{\!#1}}

\def\IE{I\mkern-6mu E}

\def\resetdisplay{\displayindent=\leftmargin
\advance\displaywidth-\leftmargin}
\makeatletter

\newbox\b@xit
\def\goodboxit#1#2#3{\def\eh{width#2}\def\ev{height#2}%
\setbox\b@xit=\hbox{\kern#1pt{#3}\kern#1pt}%
\dimen1=\ht\b@xit \advance\dimen1 by #1pt \dimen2=\dp\b@xit \advance\dimen2 by #1pt
\setbox\b@xit=\hbox{\vrule\eh height\dimen1 depth\dimen2\box\b@xit\vrule\eh}%
\setbox\b@xit=\hbox{\vtop{\vbox{\hrule\ev\box\b@xit}\hrule\ev}}%
\box\b@xit}
\def\boxit#1#2{\goodboxit{#1}{0,4pt}{#2}}
\def\displayboxit#1{\boxit4{$\displaystyle#1$}}
\def\ramollit#1#2{#1=#2#1 \mou\relax}

\def\mydisplayskips{\abovedisplayskip=.9\abovedisplayskip\mou
\belowdisplayskip=.9\belowdisplayskip\mou
\abovedisplayshortskip=.9\abovedisplayshortskip\mou
\belowdisplayshortskip=.9\belowdisplayshortskip\mou
\ramollit\baselineskip{1}
\ramollit\smallskipamount1
\ramollit\medskipamount1
\ramollit\bigskipamount1
}

\def\big#1{{\hbox{\mathsurround0pt$\left#1\vbox to10\p@{}\right.\n@space$}}}
\def\Big#1{{\hbox{\mathsurround0pt$\left#1\vbox to12\p@{}\right.\n@space$}}}
\def\bigg#1{{\hbox{\mathsurround0pt$\left#1\vbox to15\p@{}\right.\n@space$}}}
\def\Bigg#1{{\hbox{\mathsurround0pt$\left#1\vbox to17\p@{}\right.\n@space$}}}

\def\scalebigs#1{%
\def\big##1{{\hbox{\dimen0=#1pt\mathsurround0pt$\left##1\vbox to10\dimen0{}\right.\n@space$}}}%
\def\Big##1{{\hbox{\dimen0=#1pt\mathsurround0pt$\left##1\vbox to12\dimen0{}\right.\n@space$}}}%
\def\bigg##1{{\hbox{\dimen0=#1pt\mathsurround0pt$\left##1\vbox to15\dimen0{}\right.\n@space$}}}%
\def\Bigg##1{{\hbox{\dimen0=#1pt\mathsurround0pt$\left##1\vbox to17\dimen0{}\right.\n@space$}}}%
\ignorespaces}

\let\normalsizeorg\normalsize
\def\normalsize{\normalsizeorg\mydisplayskips
}

\let\smallsizeorg\small
\def\small{\smallsizeorg
\mydisplayskips
}

\def\mysmash{\relax 
  \ifmmode\def\next{\mathpalette\mathsm@mysh}\else\let\next\makesm@mysh
  \fi\next}
\def\makesm@mysh#1{\setbox\z@\hbox{#1}\finsm@mysh}
\def\mathsm@mysh#1#2{\setbox\z@\hbox{$\m@th#1{#2}$}\finsm@mysh}
\def\finsm@mysh{\ht\z@\z@ \dp\z@\z@ \box\z@ }
\def\finsm@myshbot{\dp\z@\z@ \box\z@}
\def\finsm@myshtop{\ht\z@\z@ \box\z@}
\def\finhalfsm@myshbot{\dp\z@0.5\dp\z@ \box\z@}
\def\finhalfsm@myshtop{\ht\z@0.5\ht\z@ \box\z@}
\def\finhalfsm@mysh{\ht\z@0.5\ht\z@ \dp\z@0.5\dp\z@ \box\z@}
\def\fingoodsm@mysh#1#2{\ht\z@#1\ht\z@ \dp\z@#2\dp\z@ \box\z@}
\def\smashtop#1{\begingroup\let\finsm@mysh\finsm@myshtop\mysmash{#1}\endgroup}
\def\smashbot#1{\begingroup\let\finsm@mysh\finsm@myshbot\mysmash{#1}\endgroup}
\def\halfsmashtop#1{\begingroup\let\finsm@mysh\finhalfsm@myshtop\mysmash{#1}\endgroup}
\def\halfsmashbot#1{\begingroup\let\finsm@mysh\finhalfsm@myshbot\mysmash{#1}\endgroup}
\def\halfsmash#1{\begingroup\let\finsm@mysh\finhalfsm@mysh\mysmash{#1}\endgroup}
\def\goodsmash#1#2#3{\begingroup\def\finsm@mysh{\fingoodsm@mysh{#1}{#2}}\mysmash{#3}\endgroup}

\def\mrlap#1{\rlap{$\mathsurround0pt{}#1{}$}}

\def\fixepreskip{\abovedisplayshortskip\abovedisplayskip}
\def\fixepostskip{\belowdisplayshortskip\belowdisplayskip}
\def\fixedskips{\belowdisplayshortskip\belowdisplayskip
\abovedisplayskip\belowdisplayskip\abovedisplayshortskip\abovedisplayskip}
\def\preskip{\afterassignment\fixepreskip\abovedisplayskip}
\def\postskip{\afterassignment\fixepostskip\belowdisplayskip}
\def\dskips{\afterassignment\fixedskips\belowdisplayskip}

\newcount\dimcount
\newif\ifhalfds
\def\scaledimen#1#2/#3 {\dimcount=#1\relax
\divide\dimcount by#3
\multiply\dimcount by #2
#1=\dimcount sp }

\def\scaledisplayskips#1/#2 {\scaledimen\abovedisplayskip#1/#2 
\scaledimen\abovedisplayshortskip#1/#2 
\scaledimen\belowdisplayskip#1/#2 
\scaledimen\belowdisplayshortskip#1/#2 
}

\def\halfdisplayskips{\ifhalfds\else
\scaledimen\abovedisplayskip50/100
\scaledimen\abovedisplayshortskip50/100
\scaledimen\belowdisplayskip50/100
\scaledimen\belowdisplayshortskip50/100
\halfdstrue\fi}

\def\rigiddisplayskips{\relax
\abovedisplayskip1\abovedisplayskip
\abovedisplayshortskip1\abovedisplayshortskip
\belowdisplayskip1\belowdisplayskip
\belowdisplayshortskip1\belowdisplayshortskip}

\def\mathrigid#1 {\thinmuskip=#1
\medmuskip=#1
\thickmuskip=#1 }

\thinmuskip=1mu

\def\mathemule#1#2{{\setbox0=\hbox{$\m@th#1$}\hbox to\wd0{$\m@th#2$\hss}}}
\def\mathemulerlap#1#2{{\setbox0=\hbox{$\m@th#1$}\hbox to\wd0{$\m@th#2$\hss}}}
\def\mathemulellap#1#2{{\setbox0=\hbox{$\m@th#1$}\hbox to\wd0{\hss$\m@th#2$}}}
\def\mathemulec#1#2{{\setbox0=\hbox{$\m@th#1$}\hbox to\wd0{\hss$\m@th#2$\hss}}}

\def\shift#1 #2/{\hskip#1#2\hskip-#1}
\def\kshift#1#2{\kern#1#2\kern-#1}
\def\scriptstack#1{{\let\Hbox\hbox
\setbox0=\vtop{\offinterlineskip\stackstyle
\lineskiplimit1pt\lineskip1pt
\spoolscriptstack#1\\\end\\}\relax
\edef\Hbox##1{\hbox to\wd0{\stacklhss##1\stackrhss}}\relax
\vtop{\offinterlineskip\stackstyle
\lineskiplimit1.2pt\lineskip1.2pt
\spoolscriptstack#1\\\end\\}\relax
}}

\def\stackinsidemacro#1{\scriptsize\mathsurround0pt$#1$}

\def\spoolscriptstack#1\\{\def\donne{#1}\ifx\donne\myend\let\nextspoolscriptstack
\relax\else\let\nextspoolscriptstack\spoolscriptstack
\Hbox{\stackinsidemacro{\ignorespaces#1\unskip}}\fi
\nextspoolscriptstack}

\def\cstack#1{\ifmmode\mathrel{\vcenter{\scriptsize
\tackstyle\dimen0=0pt
\def\hboxto{\setbox0=\hbox}\spoolcstack#1\\\end\\
\def\hboxto##1{\hbox to\dimen0{\stacklhss##1\stackrhss}}\spoolcstack
#1\\\end\\}}\else$\cstack{#1}$\fi}

\def\stackstyle{}
\let\stacklhss\hss
\let\stackrhss\hss

\def\spoolcstack#1\\{\def\donne{#1}\ifx\donne\myend
\let\nextspoolcstack\relax\else
\hboxto{\ignorespaces#1\unskip}\ifdim\dimen0<\wd0 \dimen0=\wd0\fi
\let\nextspoolcstack\spoolcstack\fi
\nextspoolcstack}

\postdisplaypenalty1000
\def\set#1/{\{#1\}}
\def\bigset#1/{\big\{{\let\mid\bigmid#1}\big\}}
\def\Bigset#1/{\Big\{{\let\mid\Bigmid#1}\Big\}}
\def\varset#1#2/{\left\{\vcenter to#1{}\right.\!
        {\def\mid{\varmid{#1}}#2}\!\left.\vcenter to#1{}\right\}}

\def\mid{\mathrel{|}}
\def\bigmid{\mathrel{\big|}}
\def\Bigmid{\mathrel{\Big|}}
\def\varmid#1{\mathrel{\left|\vcenter to#1{}\right.}}
\def\adaptset#1/{{\setbox0=\hbox{$\m@th#1$}\dimen0=\ht0
\advance\dimen0 by\dp0\varset{\dimen0}\vcenter{\hbox{$\m@th#1$}}/}}

\let\@addpunctorg\@addpunct
\def\nopoint{\def\@addpunct##1{\global\let\@addpunct\@addpunctorg}}

\newif\ifboldmath
\let\boldmathorg\boldmath
\def\boldmath{\protect\boldmathtrue\protect\boldmathorg}

\def\virg{\raise2pt\hbox{,}\,}

\let\cont\supseteq

\def\limite{\setbox0=\vtop{\offinterlineskip\setbox0=\hbox{lim}\dimen0=\wd0
\setbox1=\hbox{$\longrightarrow$}\ifdim\dimen0<\wd1 \dimen0=\wd1 \fi
\hbox to\dimen0{\hss\box0\hss}\kern2pt\hbox
to\dimen0{\hss\box1\hss}}\dp0=5pt\box0}
\def\limitegauche{\setbox0=\vtop{\offinterlineskip\setbox0=\hbox{lim}\dimen0=\wd0
\setbox1=\hbox{$\longleftarrow$}\ifdim\dimen0<\wd1 \dimen0=\wd1 \fi
\hbox to\dimen0{\hss\box0\hss}\kern2pt\hbox
to\dimen0{\hss\box1\hss}}\dp0=5pt\box0}

\def\confer{{\em cf.\/} \ignorespaces}
\def\idest{{\em i.e.\/} \ignorespaces}

\def\latin#1#2#3#4{\ifnum#4 = 0 \gdef#1{{\em #3 }\gdef#1{{\em #3 }}}\else
\gdef#1{{\em #2}\gdef#1{{\em #3}}}\fi}
\latin{\cf}{confer}{cf.}{0}
\latin{\Cf}{Confer}{Cf.}{0}
\latin{\confer}{confer}{cf.}{0}
\latin{\Confer}{Confer}{Cf.}{0}
\latin{\etc}{et c\char26tera}{etc}{0}
\latin{\ie}{id est}{i.e.}{0}
\latin{\Idest}{Id est}{I.e.}{0}
\latin{\Ie}{Id est}{I.e.}{0}
\parskip0.5ex
\parindent1em
\def\list#1#2{%
  \ifnum \@listdepth >5\relax
    \@toodeep
  \else
    \global\advance\@listdepth\@ne
  \fi
  \rightmargin\z@
  \listparindent1em
  \itemindent\z@
  \csname @list\romannumeral\the\@listdepth\endcsname
  \def\@itemlabel{#1}%
  \let\makelabel\@mklab
  \@nmbrlistfalse
  #2\relax
  \@trivlist
  \parskip\parsep
  \parindent\listparindent
  \advance\linewidth -\rightmargin
  \advance\linewidth -\leftmargin
  \advance\@totalleftmargin \leftmargin
  \parshape \@ne \@totalleftmargin \linewidth
  \ignorespaces}
  
\def\rien{}
\def\refstepcounter#1{\stepcounter{#1}%
\edef\previouslabel{\csname p@#1\endcsname}\ifx\previouslabel\rien
\protected@edef\@currentlabel
       {\csname the#1\endcsname}%
\else
\protected@edef\@currentlabel
       {\csname p@#1\endcsname-\csname the#1\endcsname}%
\fi
}

\def\mylistskips{\topsep4pt\relax 
         \partopsep=1\parskip \relax
         \itemsep4pt \relax
         \parsep2pt \mou\relax
         \listparindent1em
         \itemindent0pt 
         \parskip2pt\mou
         \@beginparpenalty=0
\varlistskips}

\def\nobreakitem{\@beginparpenalty=200}
\def\allowbreakitem{\@beginparpenalty=-51}

\def\mynobreak{\@nobreaktrue\@beginparpenalty=10000 }

\let\varlistskips\relax

\def\mymakelabel{\def\makelabel##1{\hss\llap{##1}}}
\def\advancelistdepth#1{\advance\@listdepth by#1}

\def\enumerate{%
  \ifnum \@enumdepth >\thr@@\@toodeep\else
    \advance\@enumdepth\@ne
    \edef\@enumctr{enum\romannumeral\the\@enumdepth}%
      \expandafter
      \list
        \csname label\@enumctr\endcsname
        {\usecounter\@enumctr\def\makelabel##1{\hss\llap{##1}}\mylistskips}%
  \fi}

\def\endenumerate{\endlist\global\let\varlistskips\relax}

\newenvironment{myenumerate}{\advance\@enumdepth\@ne
\edef\@enumctr{enum\romannumeral\the\@enumdepth}\expandafter
\list\csname label\@enumctr\endcsname
        {\usecounter\@enumctr\mymakelabel\mylistskips}}{\endlist}

\newenvironment{liste-var}[1]
{\begin{list}{#1}{\mylistskips}}{\end{list}}
\newenvironment{liste-bullet}
{\begin{list}{$\bullet$}{\mylistskips}}{\end{list}}
\newenvironment{liste-dash}
{\begin{list}{--}{\mylistskips}}{\end{list}}
\newenvironment{liste-roman}
{\begin{myenumerate}}{\end{myenumerate}}
\newenvironment{liste-Roman}
{\begin{myenumerate}}{\end{myenumerate}}
\newenvironment{liste-alpha}
{\begin{myenumerate}}{\end{myenumerate}}
\newenvironment{liste-Alpha}
{\begin{myenumerate}}{\end{myenumerate}}

\def\gob#1{}

\def\0{{\bf0}}
\def\1{{\bf1}}
\let\fonct\rightsquigarrow

\def\expression#1{``\emph{#1\/}''}

\def\itemize{%
  \ifnum \@itemdepth >\thr@@\@toodeep\else
    \advance\@itemdepth\@ne
    \edef\@itemitem{labelitem\romannumeral\the\@itemdepth}%
    \expandafter
    \list
      \csname\@itemitem\endcsname
      {\def\makelabel##1{\hss\llap{##1}}%
\labelwidth1.2em
\leftmargin\parindent\advance\leftmargin\labelwidth
\labelsep0.4em
\topsep4pt\mou
\partopsep=1ex \mou
\itemsep4pt \mou
\parsep2pt \mou\relax
         \itemindent0pt 
         \parskip0pt
         \mylistskips}\fi}

\def\enditemize{\endlist\global\let\varlistskips\relax}

\renewcommand\theenumi{\alph{enumi}}

\renewcommand\theenumii{\romannumeral\c@enumii}

\setlength\leftmargini  {1.5em}
\setlength\leftmarginii  {1.5em}
\setlength\leftmarginiii {2.5em}
\setlength\leftmarginiv  {3em}

\def\@listi{\leftmargin\leftmargini
            \parsep \parskip\relax
            \topsep 4\p@ \mou\relax
            \itemsep4\p@ \mou}
\def\@listii{\leftmargin\leftmarginii
            \parsep \parskip\relax
            \topsep 4\p@ \mou\relax
            \itemsep4\p@ \mou}
\let\@listI\@listi

\def\dobinom#1#2{\mathrel{\hbox{\small$\Big(\mkern1mu\vcenter{\footnotesize\let\scriptstyle\relax\setbox0=\hbox{$\scriptstyle#1$}
\setbox1=\hbox{$\scriptstyle#2$}
\dimen0=\wd0
\ifdim\dimen0<\wd1\dimen0=\wd1\fi
\hbox to\dimen0{\hss$\scriptstyle#1$\hss}
\hbox to\dimen0{\hss$\scriptstyle#2$\hss}
}\mkern1mu\Big)$}}}

\let\binomorg\binom
\def\binom#1#2{\mathchoice
{\dobinom{#1}{#2}}
{\binomorg{#1}{#2}}
{\binomorg{#1}{#2}}{\binomorg{#1}{#2}}}

\long\def\mysettings{
\thinmuskip=3mu
\medmuskip=4mu plus 2mu minus 2mu
\thickmuskip=5mu plus 2mu minus 2mu
\def\thinspace{\kern .16667em }
\def\negthinspace{\kern-.16667em }
\def\enspace{\kern.5em }
\def\enskip{\hskip.5em\relax}
\def\quad{\hskip1em\relax}
\def\qquad{\hskip2em\relax}
\def\smallskip{\vskip\smallskipamount}
\def\medskip{\vskip\medskipamount}
\def\bigskip{\vskip\bigskipamount}
\pretolerance=50		
\tolerance=500			
\brokenpenalty=100 			
\predisplaypenalty=10000 		
\postdisplaypenalty=7500 	
\interlinepenalty=0 		
\floatingpenalty=0 		
\hbadness=1000 				
\vbadness=1000 				
\linepenalty=10 			
\hyphenpenalty=50 			
\exhyphenpenalty=50 		
\binoppenalty=700 			
\relpenalty=500 			
\clubpenalty=5000			
\widowpenalty=5000 			
\displaywidowpenalty=3500	
\doublehyphendemerits=1000 		
\finalhyphendemerits=500 		
\scaledisplayskips80/100
\mydisplayskips
\parindent1em 
\frenchspacing
\parskip1ex plus 0.2pt minus1pt
}
\mysettings

\footnotesep=1.5ex
\def\footnoterule{\kern-3\p@
\vrule height7pt width0pt  \hrule \@width 2in \kern 2.6\p@} 

\def \newpage {%
  \if@noskipsec
    \ifx \@nodocument\relax
      \leavevmode
      \global \@noskipsecfalse
    \fi
  \fi
  \if@inlabel
    \leavevmode
  \fi
  \if@nobreak \@nobreakfalse \everypar{}\fi
  \par
  \vfil
  \penalty -\@M}

\def\centrer{
 \oddsidemargin-1truein
\evensidemargin-1truein
\dimen0=\paperwidth
\advance\dimen0 by -\textwidth
\advance\oddsidemargin by 0.5\dimen0
\advance\evensidemargin by 0.5\dimen0
}

\centrer

\def\@makefnmark{\hbox{\rm(\@textsuperscript{\normalfont\@thefnmark})\,}}
\def\@makefnmark{\hbox{$\scriptspace=0pt^{\@thefnmark}$}}
\let\puttextsuperscript\@textsuperscript

\def\signature{\vtop to0cm{\footnotesize\kern0pt\offinterlineskip%
\setbox0=\hbox{\hss ALBERTO ARABIA\hss}%
\kern5pt\hbox to\wd0{\hss ALBERTO ARABIA\hss}
\kern2mm
\hbox to\wd0{\hss CNRS--IMJ\hss}\kern2mm
\hbox to\wd0{\hss ThÈorie des Groupes\hss}
\kern2mm
\hbox to\wd0{\hss \hui\hss}
\vss}}

\def\hui{\number\day \ifnum\day=1 $^{\mathord{\rm er}}$\fi \space
\ifcase\month\or%
            janvier\or f\'evrier\or mars\or avril\or mai\or juin\or%
            juillet\or ao\^ut\or septembre\or octobre\or novembre\or %
            d\'ecembre\fi\space\number\year}

\def\makeglossary{%
  \newwrite\@glossaryfile
  \immediate\openout\@glossaryfile=\jobname.glo
  \def\glossary{\@bsphack\begingroup
                \@sanitize
                \@wrCMglossary}\typeout
    {Writing glossary file \jobname.glo }%
  \let\makeglossary\@empty
}

\def\@wrCMglossary#1{%
   \protected@write\@glossaryfile{}%
      {\string\indexentry{#1\string\fillcomma}{\thepage}}%
 \endgroup
 \@esphack}
 \def\fillcomma#1{#1\null\nobreak\hfill\nobreak}

\def\@makesindexhead#1{%
  {\parindent \z@ \raggedright
    \normalfont
    \interlinepenalty\@M
\section*{#1}
  }}
\newif\ifStyleIndexDeuxColonnes
\StyleIndexDeuxColonnestrue

\def\insideindexmacros{%
\def\ii##1{{\it##1}}%
\def\bb##1{{\bf##1}}%
\def\ul##1{{\vtop{\hbox{##1}\kern1pt\hrule}}}%
\def\ulbb##1{{\vtop{\hbox{\bf##1}\kern1pt\hrule}}}%
\def\ulii##1{{\vtop{\hbox{\it##1}\kern1pt\hrule}}}%
\raggedright}

\renewenvironment{theindex}
               {\@makesindexhead{\indexname}
                \addcontentsline{toc}{section}{\indexname}
                \thispagestyle{plain}\parindent\z@
                \parskip\z@ \@plus .3\p@\relax
                \columnseprule 0pt\small
                \spaceskip0.5ex plus1ex minus0.2ex
                \xspaceskip0.5ex plus1ex minus0.2ex
                \interlinepenalty1000
                \insideindexmacros
                \def\item{\par\penalty-50\@idxitem}}
              {\par\vskip1em\ifStyleIndexDeuxColonnes
              \if@restonecol\onecolumn\else\clearpage\fi\else\fi}
\begingroup
\catcode`\:=12
\gdef\printglossary{\immediate\closeout\@glossaryfile\begingroup\catcode`\:=12 \catcode`\@=11
\def\insideindexmacros{\lineskip2.7pt\mou\lineskiplimit2.7pt\rightskip2em} 
\def\@idxitem{\par\hangindent 30\p@}
\def\@idxitem##1:{\par\hangindent 2cm\noindent
           \setbox110=\hbox{##1\unskip.}\ifdim\wd110<2cm\hbox to2cm{\copy110\hss}\else\box110\ \fi
           \ignorespaces}
\def\subitem{\@idxitem \hspace*{2.5cm}}
\def\subsubitem{\@idxitem \hspace*{3cm}}
\def\indexspace{\par \vskip 10\p@ \@plus5\p@ \@minus3\p@\relax}
\def\indexname{Glossary}\StyleIndexDeuxColonnesfalse
\def\indexentry##1##2{\cleanalias{##1}\expandafter\item\clean\nobreak\leaders\hbox{.}\hfill\nobreak\hskip-1em\rlap{\hbox to2em{\hss##2}}}
\def\fillcomma##1{\relax}%
\parskip0pt
\begin{theindex}\input \jobname.glo\end{theindex}\endgroup
}

\gdef\cleanalias#1{\def\donne{#1}\docleanalias#1@//\ifx\donne\clean\else 
\expandafter\goblastAT\aliasresidu//\fi}
\gdef\goblastAT#1@//{\def\clean{#1}}
\gdef\docleanalias#1@#2//{\def\clean{#1}\def\aliasresidu{#2}}
\endgroup

\def\killline{\vskip-0.5\baselineskip\vskip0pt}

\long\def\@footnotetext#1{\insert\footins{%
    \reset@font\footnotesize
    \interlinepenalty\interfootnotelinepenalty
    \splittopskip\footnotesep
    \splitmaxdepth \dp\strutbox \floatingpenalty \@MM
    \hsize\columnwidth \@parboxrestore
    \protected@edef\@currentlabel{%
       \csname p@footnote\endcsname\@thefnmark
    }%
    \color@begingroup
     \@makefntext{\,%
        \rule\z@ \footnotesep 
        \ignorespaces#1\@finalstrut\strutbox}
    \color@endgroup}}
\makeatother

\def\underlinee#1{\,\vtop{\setbox0=\hbox{$\mathsurround-0,5pt#1$}\dp0=0pt\copy0\kern1.5pt\hrule\kern1pt\hrule}\,}
\def\uunderlinee#1#2{\,\vtop{\setbox0=\hbox{$\mathsurround-0,5pt#2$}\dp0=0pt\copy0\kern#1\hrule height1pt
\kern-#1\kern-1pt}\,}

\def\udp#1{\mathchoice
{\uunderlinee{1.5pt}{\displaystyle#1}}
{\uunderlinee{1.5pt}{\textstyle#1}}
{\uunderlinee{1.pt}{\scriptstyle#1}}
{\uunderlinee{1.pt}{\scriptscriptstyle#1}}
}

\def\mod{\; \mathop{\rm mod}}

\def\tiretfragile{\ifmmode\mathchoice{\hbox{-}}{\hbox{-}}{\hbox{\scriptsize-}}{\hbox{\tiny-}}\else-\fi}
\def\tiret{\protect\tiretfragile}

\def\putat#1#2#3{\vtop to0pt{\kern-#2\rlap{\kern#1#3}\vss}}
\def\[#1]{[\mkern-3mu[#1]\mkern-3mu]}

\def\partichar{{\mkern1mu\geq\mkern1mu}}

\def\parti(#1){{(
\doparti#1,,)}}
\def\doparti#1,#2,{\def\donne{#2}#1\ifx\donne\rien
\let\nextdoparti\relax
\else\partichar
\def\nextdoparti{\doparti#2,}\fi\nextdoparti}

\def\uline#1{\vtop to0pt{\mathsurround0pt\hbox{$\scriptstyle#1$}
\kern-0\prevdepth\kern1.5pt\hrule height1pt\vss}}
\def\oline#1{\vbox to0pt{\vss\mathsurround0pt\hrule height1pt\kern1.5pt
\hbox{$\scriptstyle#1$}}}

\def\dofd#1#2#3{{\muskip0=#1
\mkern-\muskip0\uline{\mkern\muskip0{#3}}#2}}

\def\usp#1{^{\dofd{1mu}{}{#1}}}

\usepackage{makeidx}\makeindex\makeglossary

\newcount\indexYNG
\newcount\rowindexYNG
\newcount\boxindexYNG
\newdimen\dimcellx \dimcellx1.5ex
\newdimen\dimcelly \dimcelly1.5ex
\newdimen\epaisseur
\newtoks\greyrows
\definecolor{mycolor}{gray}{0.70}
\greyrows={}
\newif\ifgreyYNG

\def\dimcell{\afterassignment\setdimcellx\dimcelly}
\def\setdimcellx{\dimcellx=\dimcelly}

\def\checkifgreyYNG{\greyYNGfalse\def\fin{\end}\relax
\edef\currentrow{\the\rowindexYNG}\expandafter\docheck\the\greyrows,\end,}

\def\docheck#1,{\def\donnee{#1}\ifx
\donnee\currentrow\greyYNGtrue\def\nextdocheck##1\end,{\relax}\else\ifx\donnee\fin\let\nextdocheck\relax\else\let\nextdocheck\docheck\fi\fi\nextdocheck}

\def\defcellmacro#1{\expandafter\def\csname cell#1\endcsname}

\def\defcellmacrowithin#1#2#3{\count101=#1
\advance\count101 by -1
\loop
\ifnum\count101<#2
\advance\count101 by 1
\expandafter\defcellmacro
\expandafter{\the\count101}{#3}
\repeat}

\let\everyboxYNG\relax
\let\everycellYNG\relax
\def\everyboxYNG{\expandafter\csname box\the\boxindexYNG\endcsname}
\def\everycellYNG{\expandafter\csname cell\the\boxindexYNG\endcsname}

\def\cell#1#2{\setbox2=\hbox{{\global\advance\boxindexYNG by 1\relax
\everycellYNG\ifgreyYNG\color{mycolor}\else\color{white}\fi\vrule height#2 depth0pt width#1}\llap{\vbox to#2{\vss\hbox to#1{\hss\everyboxYNG\hss}\vss}}}\relax
\setbox2=\hbox{\vrule width\epaisseur\box2\vrule width\epaisseur}%
\setbox2=\vbox{\offinterlineskip\hrule height\epaisseur \box2\hrule height\epaisseur}%
\box2}

\def\cellule{\cell\dimcellx\dimcelly}

\def\ligne#1{\hbox{\indexYNG=#1\relax 
\global\advance\rowindexYNG1\checkifgreyYNG
\ifnum\indexYNG=0\color{white}\indexYNG=1\fi\loop
\ifnum\indexYNG > 0 \cellule\kern-\epaisseur\relax \advance\indexYNG-1 \repeat}}

\def\colonne#1{\vtop{\offinterlineskip\indexYNG=#1\relax \loop
\ifnum\indexYNG > 0 \hbox{\cellule}\kern-\epaisseur \advance\indexYNG-1 \repeat}}

\def\HYOUNG#1{\def\outyoung{}{\let\ligne0\def\fin{\end}\hspool#1,\end,}%
\def\do##1{{\dimcellx##1\dimcellx \epaisseur0.4pt\,\vcenter{\offinterlineskip\outyoung}\,}}%
\mathchoice{\do{2}}{\do{1}}{\do{0.7}}{\do{0.5}}}

\def\hyoung#1{\def\outyoung{}\relax
\global\rowindexYNG=0\relax
\global\boxindexYNG=0\relax
{\let\ligne0\def\fin{\end}\hspool#1,\end,}%
\epaisseur0.4pt\,\vcenter{\offinterlineskip\outyoung}\,}

\def\hspool#1,{\def\in{#1}\ifx\in\fin\let\next\relax\else
\xdef\outyoung{\outyoung\kern-\epaisseur\ligne{#1}}%
\let\next\hspool\fi%
\next}

\def\VYOUNG#1{\def\outyoung{}{\let\colonne0\def\fin{\end}\vspool#1,\end,}%
\def\do##1{\dimcell##1\dimcell \epaisseur0.4pt\,\vcenter{\hbox{\outyoung}}\,}%
\mathchoice{\do{2}}{\do{1}}{\do{0.7}}{\do{0.5}}}

\def\vyoung#1{\def\outyoung{}{\let\colonne0\def\fin{\end}\vspool#1,\end,}%
\epaisseur0.4pt\,\vcenter{\hbox{\outyoung}}\,}

\def\vspool#1,{\def\in{#1}\ifx\in\fin\let\next\relax\else
\xdef\outyoung{\outyoung\kern-\epaisseur\colonne{#1}}%
\let\next\vspool\fi%
\next}

\newdimen\refhtr
\setbox1=\hbox{x}\refhtr\ht1

\def\down#1#2{\vtop{\kern-\refhtr\vss\hbox{$#1{#2}$}}}
\def\DOWN#1#2{{\def\do##1{\vtop{\kern-\refhtr\vss\hbox{$##1#1{#2}$}}}\relax
\mathchoice{\do\displaystyle}{\do\textstyle}{\do\scriptstyle}{\do\scriptscriptstyle}}}

\def\t{\hyoung{5,4,6,1,2,9}}
\def\v{\vyoung{5,4,6,1,2,9}}

\newcount\iii \iii=1

\loop
{\immediate\write16{-----------\the\iii-----------}
\dimcell=0,5\dimcell
\dimcell=\iii\dimcell
\noindent \llap{\the\iii]} Exemple en taille \the\dimcell $\t^{\t^{\v}}xxx\v_{\v}$ $$x\t^{\t}\v_{\v} x$$ \endgraf}
\ifnum\iii <2 \advance\iii1 \repeat

$$\dimcell1cm\vyoung{5,4,6}$$

\end

\def\Vecf{{\mathop{\bf Vec}\nolimits_f}}
\def\FI{{\mathop{\rm FI}}}
\def\FB{{\mathop{\rm FB}}}
\def\FIg{{\mathop{\rm FI}}}
\def\FBg{{\mathop{\rm FB}}}
\def\FImod{\mathop{\rm FI\tiret\rm mod}\nolimits}
\def\FBmod{\mathop{\rm FB\tiret\rm mod}\nolimits}
\def\RS{{(\mathop{\rm RS}){}}}
\def\PC{{(\mathop{\rm PC}){}}}
\def\trivial_#1{\QQ(#1)}
\def\trivial{\QQ}
\def\alternating{\epsilon\QQ}

\def\rankRS{\mathop{\rm rank}\nolimits_{\mathsc{rs}}}
\def\rankPC{\mathop{\rm rank}\nolimits_{\mathsc{pc}}}

\def\rankFIRS{\mathop{\rm rank}\nolimits^{\mathsc{fi}}_{\mathsc{rs}}}

\def\iii[#1]{[\mkern-2.75mu[#1]\mkern-2.75mu]}
\let\chiorg\chi
\def\chi{\raise2pt\hbox{$\chiorg$}} 
\def\ulambda{{\underline\lambda}}
\def\unu{{\underline\nu}}

\def\uunderline#1{\vtop{\offinterlineskip\hbox{#1}\kern1pt\hrule height1pt}}
\def\varuunderline#1#2{\vbox{\kern-0pt\offinterlineskip\hbox{#2}\kern1pt\hrule height1pt}}
\def\prodnl{\prod\nolimits}
\def\sumnl{\sum\nolimits}
\def\bigoplusnl{\bigoplus\nolimits}
\def\defmacro#1{\expandafter\def\csname#1\endcsname}

\def\arinto{\ar@{>->}}
\def\ardans{\ar@{^{(}->}}
\def\ind{\mathop{\rm ind}\nolimits}
\def\mystrut{\vrule height7pt depth4pt width0pt}
\def\weight{\mathop{\doweight}{}}
\def\doweight{
\mathchoice
{\hbox{\bfseries\itshape w}}
{\hbox{\bfseries\itshape w}}
{\hbox{\scriptsize\bfseries\itshape w}}
{\hbox{\tiny\bfseries\itshape w}}}

\let\into\rightarrowtail
\let\hook\hookrightarrow

\def\degw{\mathop{\fs{\mathop{\rm deg}}}\nolimits}
\def\degw{\mathop{{\mathop{\rm deg}\nolimits_{\weight}}}\nolimits}

\def\totdeg{\mathop{\rm deg}_{\rm tot}}
\def\tr{\mathop{\rm tr}\noindent}
\def\class{{\rm cl}}
\def\chiE_#1^#2{\chi_{\IE_{#1}^{#2}}}
\def\sqtimes{\mathbin{\mathchoice
{\scriptstyle\boxtimes}
{\scriptstyle\boxtimes}
{\scriptscriptstyle\boxtimes}
{\scriptscriptstyle\boxtimes}
}}
\def\siminfty{\mathbin{\vtop{\offinterlineskip\setbox0=\hbox{$\simeq$}
\copy0\kern2pt
\hbox to \wd0{\hss$\scriptscriptstyle\infty$\hss}
}}}
\let\iGamma\varGamma
\def\mod#1{{\, (\mathop{\rm mod}#1)}}

\def\supp|#1|{\{\mkern-6mu\{#1\}\mkern-6mu\}}
\def\itilde#1{\,\nobreak\tilde{\!#1}}
\def\mycases#1{\left\{\mathalign{#1}\right.}

\def\varscalar#1#2#3#4{#1\langle\mkern-#2mu#1\langle{\,}#3\mathbin{#1|}#4{\,}#1\rangle\mkern-#2mu#1\rangle}
\def\scalar{\varscalar{}{4}}
\def\bigscalar{\varscalar{\big}{6}}

\def\inhibe{\begingroup\color{red}}
\let\endinhibe\endgroup
\long\def\inhibe#1\endinhibe{}

\def\tthref#1#2{{\tt#1}}
\def\varhref{\ifx\href\undefined
\let\nexthref\tthref\else
\let\nexthref\href\fi\nexthref}
\ifarxiv
\else
\usepackage{hyperref}
\hypersetup{
    colorlinks=true,
    linkcolor=blue,
    filecolor=magenta,      
    urlcolor=cyan,
}
\fi
\title{On the equivalence of two stability conditions\\ of $\FB$-modules}
\author{Alberto Arabia\footnote{UniversitÈ Paris-Diderot, IMJ-PRG, CNRS, B‚timent Sophie Germain, bureau 608, Case 7012, 75205. Paris Cedex 13, France. Contact\,: {\tt alberto.arabia@imj-prg.fr}.}}
\date{May 16, 2018}
\begin{document}
\mysettings
\englishtrue

{\parskip0pt\maketitle}


\kern-\parskip\vskip-0em
{\small\noindent
{\bf Abstract.} 
\parskip0.25ex \mou
\parindent0pt
\rigiddisplayskips
\scaledisplayskips60/100
We give a proof of the equivalence of two stability conditions in the work of Thomas Church and Benson Farb  on $\FB$-modules (\cite{chu-far})\,:  \emph{representation stability}~$\RS$ and \emph{character polynomiality} $\PC$. While the implication $\RS \Rightarrow \PC$ is a simple consequence of the Frobenius character formula, the converse seems not to have been documented, which motivated this work.

\smallskip
An $\FB$-module is a countable family $\W:=\set W_m/_{m\in\NN}$ 
of linear representations $W_m$ of the symmetric groups $\SSS_m$, assumed, in these notes, to be finite dimensional over the field of rational numbers~$\QQ$, hereafter called \emph{an $\SSS_{m}$-module\/}\glossary{$\SSS_m$-module:a finite dimensional linear representation of the symmetric group $\SSS_m$}. 


\smallskip
$\bullet$ The $\FB$-module $\W$ is said to be \emph{representation stable} (RS in short), if there exists some $N\in\NN$, such that, for all $m\geq N$,
the decomposition in simple factors of $W_m$ has the \emph{stable} form\,:
$$
\postskip-0.5ex
\hskip1.55cm
W_{m}\sim\bigoplusnl_{\lambda} V_{\lambda[m]}{}^{\ng_\lambda}\,,$$
where 
\begin{list}{\raise1pt\hbox{$\scriptstyle\rhd$}}{\topsep4pt\parsep0pt\itemsep1pt}
\item 
$\lambda:=\parti(\lambda_1,\ldots,\lambda_\ell)$ is a partition verifying $|\lambda|+\lambda_1\leq N$; 
\item $V_{\lambda[m]}$ is the simple $\SSS_{m}$-module associated with 
$\mathrigid3mu
\lambda[m]:=\parti(m{-}|\lambda|,\lambda_1,\ldots,\lambda_\ell)$;
\item  for $m\geq N$, the family of natural numbers
$\set \ng_{\lambda}/_{\lambda}$ is \emph{constant}.
\end{list}

The smallest such $N$ is \emph{the rank of $\RS$ of $\W$}, which we will denote by 
`$\rankRS(\W)$'.

\smallskip
$\bullet$ The $\FB$-module $\W$ is said to \emph{have a polynomial character} (PC in short), if there exists $N\in\NN$ and a polynomial $P_{\W}\in\QQ[X_1,\ldots,X_N]$, such that, for all $m\geq N$ and all $g\in\SSS_m$, one has\,:
$$\preskip-0ex
\chi_{W_m}(g)=P_{\W}(X_1(g),\ldots,X_N(g))\,,
\postskip1.5ex
$$
where $\chi_{W_m}$ is the character of  $W_m$, and 
$X_i(g)$ is the number of cycles of length $i$ in the decomposition of $g$ in $\SSS_m$ as a product of disjoint cycles.
The polynomial~$P_\W$, which is unique, is \emph{the polynomial character of $\W$}. The smallest such $N$ is \emph{the rank of character polynomiality of $\W$}, it will be denoted by `$\rankPC(\W)$'.

\smallskip
Regarding the equivalence of these conditions, we prove the following theorem.

\medskip
\noindent{\bf Theorem (\ref{theo-PC=RS}). }{\slshape Let $\W$ be an $\FB$-module.
\def\varlistskips{\advance\leftmargin1em\topsep1pt\itemsep2pt}\begin{enumerate}
\mynobreak\nobreak
\item 
$\rankPC(\W)\leq\rankRS(\W)$.
\item If $\W$ is $\PC$ and has polynomial character $P_\W$, then 
$$
\rankRS(\W)\leq \max\set \rankPC(\W),2\degw(P_{\W})/\,,
$$
where $\degw(P_\W)$ denotes the polynomial degree of $P_\W$, when stipulating that $\degw (X_i):=i$.
\end{enumerate}}

\bigskip
\noindent {\slshape Comments. }(i) Theorem (\ref{prop-basic-CP}) exhibit $\FB$-modules $\W$ with $\rankPC(\W)=0$ and arbitrarily big  $\rankRS(\W)$. 
It also shows that the bounds in (\ref{theo-PC=RS}) are optimal. (ii)
Beyond $\FB$-modules, we also look at $\FI$-modules $\W$,  which are $\FB$-modules with additional structure (\ref{ssec-FI-modules}). There is a corresponding representation stability condition for $\FI$-modules, its rank is denoted by $\rankFIRS(\W)$ (\ref{sssec-RS-FI-modules}). One has $\rankRS(\W)\leq\rankFIRS(\W)$, but there is no equivalence with $\PC$ stability, in particular thm. \ref{theo-PC=RS}-(\ref{theo-PC=RS-b}) for $\FI$-modules and $\rankFIRS$ is false.

\smallskip
An interesting  by-product in this work is 
the fact that the \emph{weight} $\weight(W_m)$ of an $\SSS_m$-module $W_m$ (see \ref{sec-weight}) coincides with the lowest possible degree $\degw(P)$ of a polynomial $P$ expressing the character $\chi_{W_m}$  (\ref{theo-degree-poly}). Moreover, if $\W:=\set W_m/$ is a $\PC$ $\FB$-module, then 
$\halfdisplayskips\degw(P_{\W})=
\weight(W_m)\,,
$
for all $m\geq\rankRS(\W)$.
The next theorem follows from these observations.

\noendpoint\begin{prop*}[\ref{prop-weight-tensor-product}-(\ref{prop-weight-tensor-product-b}).]
If $\W_1$ and $\W_2$ are $\FB$-modules which are $\PC_{m\geq N}$, then 
$$\weight(W_{1,m}\otimes W_{2,m})=\weight(W_{1,m})+\weight(W_{2,m})\,,$$
for all $m\geq 2(\weight(W_{1,m})+\weight(W_{2,m}))$. In particular,
$$\weight((\W_{1}\otimes\W_{2})_{\geq N})=
\weight((\W_{1})_{\geq N})+
\weight((\W_{2})_{\geq N})\,,
$$
where $(\W)_{\geq N}$ denotes the \emph{truncation} of $\W$ that replaces $W_n$ by $0$ for all $n<N$.
\end{prop*}

\setcounter{tocdepth}{2}
\tableofcontents

}

\vskip1em\noindent
The first two sections are elementary. They introduce notations and well-known properties of $\FI$ and $\FB$-modules, which advanced readers can skip.

\vskip-1em\vskip0pt\section{Preliminaries}
\vskip-1ex\vskip0pt\subsection{General notations}
\def\varlistskips{\leftmargin1.em\labelwidth1cm
\labelsep3pt\topsep3pt\itemsep1pt\mou\parsep1pt}\begin{itemize}
\item Given a group $G$, we denote by $G\tiret\rm mod$ the category of 
$G$-modules, \idest of finite dimensional linear representations of $G$ over the field of rational numbers $\QQ$.
\item We denote by $\SSS_m$ the \emph{symmetric group} defined over $\iii[1,m]$. Its elements are the bijections of the interval of natural numbers $\iii[1,m]$. For $n\leq m$, we write $\SSS_n\dans\SSS_m$  for the inclusion which identifies a permutation $g$ of $\iii[1,n]$ with its obvious extension to $\iii[1,m]$ that fixes all $i>n$.
\item $\trivial_m$ and $\alternating_m$ denote respectively the \emph{trivial} and the \emph{alternating or signature} $\QQ$-linear representation of the group $\SSS_m$.
\item $\SSS_{a}\sqtimes\SSS_{b}$ denotes the subgroup of $\SSS_{a+b}$ stabilizing the subset $\iii[1,a]$, or, equivalently, the subset $\iii[a{+}1,a{+}b]$. The subgroup $\1_a\boxtimes\SSS_{b}\dans\SSS_a\boxtimes\SSS_{b}$ is the \emph{fixator} of $\iii[1,a]$, \idest the set of permutations fixing every $i\leq a$.
\item If $W_a$ is an $\SSS_a$-module and $W_b$ is an $\SSS_b$-module, we denote by $W_a\sqtimes W_b$ the tensor product $W_a\otimes_\QQ W_b$ endowed with the $\SSS_a\sqtimes\SSS_b$-module structure defined by the componentwise action $(g_a,g_b)(w_a\otimes w_b):=g_a(w_a)\otimes g_b(w_b)$.
\item A \emph{partition} of $m\in\NN\moins\set0/$ is any decreasing sequence of natural numbers $\lambda:=\parti(\lambda_1,\ldots,\lambda_\ell{\,>\,}0)$ such that $m=\sum_i\lambda_i$. It is also denoted as the $m$-tuple $(1^{\ng_1},2^{\ng_2},\ldots, m^{\ng_m})$ where $\ng_k:=\#\set i\mid \lambda_i=k/$, so that $m=\sum_i i\,\ng_i$. The notation $\lambda\vdash m$ says that $\lambda$ is a partition of $m$, and $|\lambda|$ is used for the number partitioned by $\lambda$. The partition $\lambda$ is \emph{empty\/} if $|\lambda|=0$.
\item $\QQ_{\class}(\SSS_m)$ denotes the $\QQ$-algebra of \emph{class functions} of $\SSS_m$. These are the functions $f:\SSS_{m}\to \QQ$ constant along the conjugacy classes of $\SSS_m$, \idest such that $f(gxg^{-1})=f(x)$, $\forall g,x\in\SSS_m$. The scalar product 
$$\scalar{\_}{\_}_{\SSS_m}:\QQ_{\class}(\SSS_m)\times \QQ_{\class}(\SSS_m)\to \QQ\postskip2pt$$
is defined by
$$\preskip2pt\hskip2.7cm
\scalar{f_1}{f_2}_{\SSS_m}:=\frac1{|\SSS_m|}\sum_{g\in\SSS_m}f_1(g)\,f_{2}(g^{-1})\,,\relax{\quad\forall f_1,f_2\in
\QQ_{\class}(\SSS_m)\,.}
$$
\item If $W_m$ is an $\SSS_{m}$-module, $\chi_{W_m}:\SSS_m\to\QQ$ denotes its character.
The \emph{Schur's orthogonality relations\/} state that if $V_1$ and $V_2$ are \emph{simple\/} $\SSS_m$-modules, then
$\scalar{\chi_{V_1}}{\chi_{V_2}}_{\SSS_m}$ is equal to $1$ if $V_1$ is isomorphic to $V_2$, and  to $0$ otherwise.
\end{itemize}

\medskip\noindent
Between here and the end of this preliminary section we recall concepts and terminology from the works of Church and Farb (\cf\cite{chu-far}).\killline

\subsection{The categories of $\FB$ and $\FI$ modules}\label{ssec-FI-modules}\vskip10pt

\def\varlistskips{\topsep4pt}
\begin{itemize}
\item {\boldmath$\FB$} denotes the category of \uunderline Finite sets and \uunderline Bijections. An \emph{$\FB$-module} is, by definition, a \emph{covariant} functor from the category $\FBg$ to the category $\Vecf(\QQ)$ of finite dimensional $\QQ$-vector spaces and $\QQ$-linear maps\,:
$$\W:\FBg\fonct\Vecf(\QQ)\,.$$
To give an $\FB$-module $\W$ is equivalent to giving the countable collection 
$\set W_m:=\W(\iii[1,m])/_{m\in\NN}$,
where $W_m$ is
an $\SSS_m$-module. A morphism of $\FB$-modules $f:\W\to\Z$ corresponds to a family $\set f_m:W_m\to Z_m/_{m\in\NN}$ of morphisms of $\SSS_m$-modules. We have a canonical identification\,:
$$\Mor_\FB(\W,\Z)=\prodnl_{m\in\NN}\Hom_{\SSS_m}(W_m,Z_m)$$

\noindent{\slshape Notation. }The category of $\FB$-modules will be denoted by $\FB$-mod. It is a semi-simple abelian category.

\item {\boldmath$\FI$} denotes the category of \uunderline Finite sets and \uunderline Injections. An \emph{$\FI$-module} is a \emph{covariant functor} from $\FIg$ to $\Vecf(\QQ)$\,:
$$\W:\FIg\fonct\Vecf(\QQ)\,.$$

\noindent To give an $\FI$-module is equivalent to giving
\def\varlistskips{
\parsep0pt
\itemindent0.5cm\labelwidth1cm\labelsep1ex\def\theenumi{FI-\arabic{enumi}}}\begin{enumerate}

\item an $\FB$-module $\W:=\set W_m/_{m\in\NN}$;

\item for all  $m\in\NN$, an \emph{interior} map
$\phi(\W)_m:W_m\to W_{m{+}1}$ (in short $\phi_m$), which is $\QQ$-linear and
such that, for all $g\in \SSS_m\dans \SSS_{m+1}$, one has
$$\phi_m(g\cdot w)=g\cdot\phi_m(w)\,.$$

\item for all $n\geq m$, the image of $\phi_{n,m}:=\phi_{n{-}1}\circ\cdots\circ\phi_m$ satisfies\,:
$$\phi_{n,m}(W_m)\dans (W_n)^{\1_m\sqtimes \SSS_{n{-}m}}\,.$$
\end{enumerate}

\noindent Under this equivalence, a morphism of $\FI$-modules $f:\W\to\Z$ is simply a morphism of $\FB$-modules which is compatible with the interior maps $\phi_m$, \idest such  that the diagrams
$$\halfdisplayskips
\xymatrix@R=.5cm@C=1.2cm{
W_m\ar[r]^(.45){\ \phi(\W)_m\ }\ar[d]_{f_m}&W_{m+1}\ar[d]^{f_{m+1}}\\
Z_m\ar[r]^(.45){\ \phi(\Z)_m\ }&Z_{m+1}
}$$
are commutative.

\noindent{\slshape Notation. }The category of $\FI$-modules will be denoted by $\FImod$. It is an abelian category, which is \emph{not\/} semi-simple.

\end{itemize}

\noendpoint\begin{comms}\label{comms-sub-cat}
\begin{enumerate}\itemsep0pt
\item\label{comms-sub-cat-a}The category $\FBmod$ is equivalent to the full subcategory of $\FImod$ of $\FI$-modules whose interior maps are null. The \emph{forgetful functor} associates to an $\FI$-module $\W:=\set \phi_m:W_m\to W_{m+1}/$ the $\FB$-module $\W:=\set W_m/$. It is an additive and exact functor which will be tacitly used.

\item\label{comms-sub-cat-b}More interesting, the subcategory $\FImod'$ of $\FI$-modules whose interior maps are \emph{injective} and \emph{(eventually) exhaustive}, \idest such that the image $\phi_m(W_m)$ generates $W_{m+1}$ as $\SSS_{m+1}$-module for large $m$. Among these, the $\FI$-modules $\V_\lambda$'s which bind together the simple modules $V_{\lambda[m]}$ (\ref{prop-Vlambda}). 
We will see that $\FB$-modules $\W$ which are $\PC$ or $\RS$ are asymptotically isomorphic to direct sums of $\V_\lambda$'s (\ref{last-rema}).
\end{enumerate}
\end{comms}

\subsubsection{The truncations}\label{comm-trunc}
(\footnote{This truncations are frequently called \emph{stupid} or \emph{brutal}.}) The functor $(\_)_{\geq \ell}:\FImod\fonct \FImod$ that 
\expression{truncates} an $\FI$-module $\set W_m/_{m\in\NN}$ by replacing  by $\0$ its terms $W_m$ for $m<\ell$, is an additive exact functor. The \emph{cokernel} of the natural inclusion  $(\_)_{\geq\ell}\hook\id_{\FIg}$   is the  truncation 
$(\_)_{< \ell}$, it replaces by $\0$ the terms $W_m$ for $m\geq \ell$ . We thus have short exact sequences
$$\0\to(\W)_{\geq\ell}\into\W\onto(\W)_{<\ell}\to\0\,,$$
which are natural with respect to $\W$.
The full subcategory $\FImod_{\geq \ell}$ of $\FI$-modules $\W$ such that the inclusion $\W_{\geq \ell}\into\W$ is an isomorphism is an abelian subcategory, and similarly for the full subcategory $\FImod_{< \ell}$ of $\FI$-modules $\W$ such that the quotient $\W\onto \W_{< \ell}$ is an isomorphism. One has, 
$$\Ext^{i}_{\FI}\big(\FImod_{\geq \ell},\FImod_{\leq \ell}\big)=0\,,\quad\forall i>0\,.$$
The intersection $\FImod_{\geq \ell}\cap \FImod_{\leq \ell}$ is  the (semi-simple) category $\SSS_\ell\tiret\rm mod$.

\subsection{Projective $\FI$-modules}\label{ssec-proj}An obvious way to construct $\FI$-modules of the type described in \hbox{\ref{comms-sub-cat}-(\ref{comms-sub-cat-b})} is to start with ‡ representation $W_n$ of some $\SSS_n$, and then set, for all $m\in\NN$\,:
$$
\P(W_n)_{m}:=\begin{cases}
\noalign{\kern-2pt}
\, \0\,,&\hbox{if $m<n$,}\\
\, \ind^{\SSS_m}_{\SSS_{n}\sqtimes \SSS_{m-n}}\big(W_n\sqtimes \trivial_{m-n}\big)\,,&\hbox{otherwise.}
\\\noalign{\kern-2pt}
\end{cases}$$
For each  $m\geq n$, the composition of the following natural maps $\iota_m$ and $\kappa_{m+1}$\,:
$$\preskip2pt
\xymatrix@R=5mm{
W_n\sqtimes \trivial_{m-n}\ardans[r]|(0.47){\ \iota_m\ }
\ar@{.>}[rd]+<-55pt,10pt>|(0.65){\mystrut\ \psi_m\ }&
W_n\sqtimes \trivial_{m+1-n}\strut\ardans[d]^{\ \kappa_{m+1}}\\
&\relax{\ind^{\SSS_{m+1}}_{\SSS_n\sqtimes \SSS_{m+1-n}}}\relax{\big(W_n\sqtimes \trivial_{m+1-n}\big)}\kern-10mm
}$$
gives the map $\psi_{m}$ whose image is invariant under 
$(\1_{n}\sqtimes \SSS_{m+1-n})$. In particular, $\psi_m$ is 
$(\SSS_{n}\sqtimes \SSS_{m-n})$-linear, inducing thus the map
$$\phi_m:=\ind(\psi_m):
\ind^{\SSS_{m}}_{\SSS_n\sqtimes \SSS_{m-n}}
\big(W_n\sqtimes \trivial_{m-n}\big)
\to
\ind^{\SSS_{m+1}}_{\SSS_n\sqtimes \SSS_{m+1-n}}
\big(W_n\sqtimes \trivial_{m+1-n}\big)\,,$$
which satisfies the requirements making the family
$$\P(W_m):=\set \phi_{m}:\P(W_n)_{m}\to\P(W_n)_{m+1}/_{m\in\NN}\,,$$
an $\FI$-module.

All the previous constructions are natural on the category of $\SSS_{n}$-modules. As a consequence, we have defined a functor
$$\P:\SSS_n\tiret{\rm mod}\fonct\FImod\,,$$
which is clearly additive and exact.

\begin{prop}\label{proj-FI-modules}Let $W_n$ be a representation of $\SSS_n$.
\def\varlistskips{\parsep0pt\itemsep2pt}\begin{enumerate}
\item The interior maps of $\P(W_n)$ are injective. They are exhaustive $\forall m\geq n$.
\item There is a natural identification of functors
$$\Mor_{\FIg}(\P(W_n),\_)=\Hom_{\SSS_{n}}(W_n,(\_)_n)\,.$$
\item The $\FI$-module $\P(W_n)$ is a projective $\FI$-module, which moreover is simple projective if and only if $W_n$ is a simple $\SSS_n$-module. All simple projective $\FI$-modules are of this form.
\item The category $\FImod$ has enough projective objects.
\end{enumerate}
\end{prop}
\demo Left to the reader.
\enddemo
\subsection{Young diagrams and Pieri's rule}
We  recall the re-parametrization of irreducible representations of the symmetric groups due to Church and Farb.

\subsubsection{The \emph{socle} and the \emph{weight} of a partition}\label{sssec-weight-partition}
Let $\lambda=(\lambda_1\geq \lambda_2\geq\cdots\geq\lambda_\ell>0)$ be a non-empty partition of $m:=\sum_i\lambda_i$.
\def\varlistskips{\parsep0pt\itemsep0pt}\begin{itemize} 
\item The \emph{socle\/} of $\lambda$ is the sub-partition $\ulambda:= \parti(\lambda_2,\dots,\lambda_\ell)$. 
\item The \emph{weight\/} of $\lambda$ is the number $\weight(\lambda):=|\ulambda|=\lambda_2+\cdots+\lambda_\ell$.
\end{itemize}
\noindent Notice that the map $\lambda\mapsto\ulambda$ is injective from the set of partitions of $m\in\NN$, so that it amounts the same, giving $(\lambda\vdash m)$ or giving $(m,\ulambda)$ with $m\geq |\ulambda|+\ulambda_1$. 
In terms of Young diagrams, in order to get the socle $\ulambda$ of $\lambda$, one simply erases the first row of $\lambda$. The picture is thus\,: 
$$\preskip1em
\dimcell1.75mm
\lambda:={\hyoung{7,4,3,1}}\quad\mapsto \quad\ulambda:=\hyoung{0,4,3,1}\,.
\vtop to0pt{\offinterlineskip\kern-.5pt
\smash{\llap{\vrule height 11pt depth-10.6pt width1.8cm \hskip2.85cm}}
\kern1pt
\smash{\llap{\vrule height 11pt depth-10.6pt width1.8cm \hskip2.9cm}}
\vss}
$$
Conversely, given $\lambda=\parti(\lambda_1,\lambda_2,\dots,\lambda_\ell)$ and any $m\geq|\lambda|+\lambda_1$, Church and Farb introduce the notation
$$
\displayboxit{\lambda[m]:=\parti(m-|\lambda|,\lambda_1,\ldots,\lambda_\ell)}
$$
which, in terms of Young diagrams, corresponds to adding a first row with as many boxes as is necessary to raise the total number of boxes from $|\lambda|$ to~$m$. 
$$\preskip1em
\dimcell1.75mm
\lambda:={\hyoung{0,7,4,3,1}}\quad\mapsto \quad\lambda[m]:={\greyrows{1}\hyoung{10,7,4,3,1}}\,.
$$
Notice the following obvious notational facts\,:
$$\lambda[m]\vdash m\,,\quad
\underline{\lambda[m]}=\lambda\,,\quad
\lambda=\ulambda[\,|\lambda|\,]\,.
$$

\subsubsection{Irreducible representations of the symmetric groups}\label{ssec-irr-Sm-Q}The irreducible representations of $\SSS_m$ are parametrized by the partitions $\nu\vdash\nobreak n$. The simple $\SSS_m$-module associated with $\nu$ is classically denoted by $V_\nu$.
The following proposition recalls a 
fundamental fact about the representations of the symmetric groups (see \cite{fulhar}, chap. 4, thm. 4.3, pp. 46--).

\begin{prop}\label{prop-irr-Sm-Q}The irreducible representations of the symmetric groups over a field $k$ of characteristic zero are defined over the field of rational numbers. If $W_m$ is a $k[\SSS_n]$-module of finite dimension, the values of the character 
$$\chi_{W_m}:\SSS_m\to k\,,\quad g\mapsto \tr(g:W_m\to W_m)\,,$$
belong to the ring of integers $\ZZ\dans k$.
\end{prop}
\def\demoname{\it Hint. }
\demo 
Because $W_m$ is defined over $\QQ$, the trace of $g\in\SSS_m$ belongs a priori to~$\QQ$, and because the trace is the sum of the eigenvalues, which are roots of unity, it is also an algebraic integer. It is therefore an integer.
\enddemo

\subsubsection{Pieri's rule}\label{sssec-Pieri}
Pieri's rule (\footnote{For a thorough introduction to these rules, read  
 \myS4.3, p.~54--62, and also Appendix~A on the \emph{Littlewood-Richardson rules}, p. 451, both in Fulton-Harris' book \cite{fulhar}.}) gives the irreducible factors of the terms $\P(V_{\nu})_m$, for any given partition $\nu\vdash n\leq m$. The rule says that in the decomposition
$$\P(V_\nu)_{m}:=\ind^{\SSS_m}_{\SSS_n\sqtimes \SSS_{m-n}}
\big(V_\nu\sqtimes \trivial_{m-n}\big)
=\bigoplus\nolimits_{\mu\vdash m}V_{\mu}{}^{\ng_{\nu}(\mu)}\,,\eqnnb\label{eq-pieri}$$
the nonzero multiplicities $\ng_{\nu}(\mu)$ are all equal to $1$, and the corresponding partitions $\mu$ are those obtained from $\nu$ by adding $(m{-}n)$ boxes in \emph{different\/} columns.
For example, if $\nu:=(3,2,2)$ and $m\in\set8,9,10,11/$, we have
$$\dimcell1.5mm
\mathalign{
\ind^{\SSS_{8}}_{\SSS_7\sqtimes \SSS_1}
V_{\hyoung{3,2,2}}&=&
{
\defcellmacro{4}{\greyYNGtrue}
V_{\hyoung{4,2,2}}}
\hfill&\oplus&
{
\defcellmacro{6}{\greyYNGtrue}
V_{\hyoung{3,3,2}}}
\hfill&\oplus&
{
\defcellmacro{8}{\greyYNGtrue}
V_{\hyoung{3,2,2,1}}}
\hfill\\
\ind^{\SSS_{9}}_{\SSS_7\sqtimes \SSS_2}
V_{\hyoung{3,2,2}}&=&
{\defcellmacro{5}{\greyYNGtrue}
\defcellmacro{4}{\greyYNGtrue}
V_{\hyoung{5,2,2}}}
\hfill&\oplus&
{\defcellmacro{4}{\greyYNGtrue}
\defcellmacro{7}{\greyYNGtrue}
V_{\hyoung{4,3,2}}}
\hfill&\oplus&
{\defcellmacro{4}{\greyYNGtrue}
\defcellmacro{9}{\greyYNGtrue}
V_{\hyoung{4,2,2,1}}}
\hfill&\oplus&
{\defcellmacro{6}{\greyYNGtrue}
\defcellmacro{9}{\greyYNGtrue}
V_{\hyoung{3,3,2,1}}}
\hfill&\oplus&
{\defcellmacro{8}{\greyYNGtrue}
\defcellmacro{9}{\greyYNGtrue}
V_{\hyoung{3,2,2,2}}}\hfill\\
\ind^{\SSS_{10}}_{\SSS_7\sqtimes \SSS_3}
V_{\hyoung{3,2,2}}&=&
{
\defcellmacrowithin{4}{6}{\greyYNGtrue}
V_{\hyoung{6,2,2}}}
\hfill&\oplus&
{\defcellmacro{4}{\greyYNGtrue}
\defcellmacro{5}{\greyYNGtrue}
\defcellmacro{8}{\greyYNGtrue}
V_{\hyoung{5,3,2}}}
\hfill&\oplus&
{
\defcellmacro{4}{\greyYNGtrue}
\defcellmacro{5}{\greyYNGtrue}
\defcellmacro{10}{\greyYNGtrue}
V_{\hyoung{5,2,2,1}}}
\hfill&\oplus&
{\defcellmacro{4}{\greyYNGtrue}
\defcellmacro{7}{\greyYNGtrue}
\defcellmacro{10}{\greyYNGtrue}
V_{\hyoung{4,3,2,1}}}
\hfill&\oplus&
{\defcellmacro{4}{\greyYNGtrue}
\defcellmacro{9}{\greyYNGtrue}
\defcellmacro{10}{\greyYNGtrue}
V_{\hyoung{4,2,2,2}}}
\hfill&\oplus&
{\defcellmacro{6}{\greyYNGtrue}
\defcellmacro{9}{\greyYNGtrue}
\defcellmacro{10}{\greyYNGtrue}
V_{\hyoung{3,3,2,2}}}\hfill\\
\ind^{\SSS_{11}}_{\SSS_7\sqtimes \SSS_4}
V_{\hyoung{3,2,2}}&=&
{
\defcellmacrowithin{4}{7}{\greyYNGtrue}
V_{\hyoung{7,2,2}}}
\hfill&\oplus&
{\defcellmacrowithin{4}{6}{\greyYNGtrue}
\defcellmacro{9}{\greyYNGtrue}
V_{\hyoung{6,3,2}}}
\hfill&\oplus&
{
\defcellmacrowithin{4}{6}{\greyYNGtrue}
\defcellmacro{11}{\greyYNGtrue}
V_{\hyoung{6,2,2,1}}}
\hfill&\oplus&
{\defcellmacro{4}{\greyYNGtrue}
\defcellmacro{5}{\greyYNGtrue}
\defcellmacro{8}{\greyYNGtrue}
\defcellmacro{11}{\greyYNGtrue}
V_{\hyoung{5,3,2,1}}}
\hfill&\oplus&
{
\defcellmacro{4}{\greyYNGtrue}
\defcellmacro{5}{\greyYNGtrue}
\defcellmacro{10}{\greyYNGtrue}
\defcellmacro{11}{\greyYNGtrue}
V_{\hyoung{5,2,2,2}}}
\hfill&\oplus&
{
\defcellmacro{4}{\greyYNGtrue}
\defcellmacro{7}{\greyYNGtrue}
\defcellmacro{10}{\greyYNGtrue}
\defcellmacro{11}{\greyYNGtrue}
V_{\hyoung{4,3,2,2}}}\hfill\\
}
$$
where, after a key observation of Church and Farb, for $m\geq|\nu|+\nu_{1}$(${}=10$ in this example), the set of socles of the Young diagrams appearing in the decomposition \ref{eq-pieri} is constant.
The following lemma obviously follows.\killline

\begin{lemm}\label{lemm-pieri}Let $\nu\vdash n>0$.
\def\varlistskips{\itemsep0.4cm}
\begin{enumerate}
\item
\label{lemm-pieri-a} For all $\mu\vdash m\geq n$, such that $\ng_{\nu}(\mu)\not =0$, the weight of $\mu$ verifies 
$$
|\unu|=\weight(\nu)\leq\weight(\mu)\leq|\nu|\,,$$ and \hfil
$\left\{\mathalign{
\weight(\mu)&=&\weight(\nu)&\ \Longleftrightarrow\ &\mu=\unu[m]\hfill\\
\weight(\mu)&=&|\nu|\hfill&\Longleftrightarrow&\big(m\geq|\nu|+\nu_1\big)\ \&\ \big(\mu=\nu[m]\big)\,.
}\right.$
\item\scaledisplayskips80/100
\label{lemm-pieri-b}
Let $m_0=|\nu|+\nu_1$ and let $\P_{\nu}$ be a set of partitions $\mu$ such that
$$
\kern-5pt
\ind^{\SSS_{m_0}}_{\SSS_n\sqtimes \SSS_{m_0-n}}\big(V_\nu\sqtimes \trivial_{m_0-n})\big)=
\bigoplus\nolimits_{\mu\in \P_{\nu}}V_{\mu[m_0]}\,.
$$
Then, $\nu[m_0]\in\P_{\nu}$ and
$\nu[m_0]_1=\nu[m_0]_2\,.$

\smallskip
\noindent For all $m> m_0$, we have
$$
\ind^{\SSS_{m}}_{\SSS_n\sqtimes \SSS_{m-n}}\big(V_\nu\sqtimes \trivial_{m-n})\big)=
\bigoplus\nolimits_{\mu\in \P_{\nu}}
V_{\mu[m]}\,.
$$
and, for all $\mu\in\P_{\nu}$, 
$\mu[m_0]_1>\mu[m_0]_2\,.$
\end{enumerate}
\end{lemm}
\subsection{The \emph{weight} of a representation and of an $\FB$-module}
\label{sec-weight}The definition of \emph{weight}, is extended
from partitions (\ref{sssec-weight-partition}) to representations and, more generally, to $\FB$-modules.
\begin{itemize}
\item The \emph{weight of an $\SSS_m$-module $W_m$} is the upper bound of the weights of the partitions associated with its irreducible factors, \idest if 
$$
W_m\sim\bigoplusnl_{\mu\vdash m} V_{\tau}{}^{\ng_{\mu}}\,,\postskip0pt$$
then 
$$\weight(W_m):=\sup\bigset\weight(\mu)\mid\ng_{\mu}\not=0/\,.$$

\medskip\noindent For example, as a consequence of lemma \ref{lemm-pieri}-(\ref{lemm-pieri-a}), we have
$$
\left\{\mathalign{
|\unu|=\weight(\nu)\leq\weight(\P(V_{\nu}))_m\leq|\nu|\,,&\quad\forall m\geq|\nu|\,,\text{ and }\hfill\\
\hfill\weight(\P(V_{\nu}))_m=|\nu|\,,&\quad\forall m\geq |\nu|+\nu_1\,.\hfill}\right.
\eqno(\ast)$$

\item The \emph{weight of an $\FB$-module $\W:=(W_m)$} is the upper-bound of the weights of its terms, \idest
$$\weight(\W):=\sup\bigset \weight(W_m)/_{m\in\NN}\,.
$$

\item The \emph{weight at infinity of an $\FB$-module $\W:=(W_m)$} is
$$\weight_{\infty}(\W):=\lim_{N\mapsto+\infty}\weight(\W_{\geq N})\,.
$$

\noindent --It is easy to see that $\weight_{\infty}(\W)=\weight(\W_{\geq N})$, for some $N\gg0$.

\noindent--We have
$
\weight(\P(V_{\lambda}))=\weight_{\infty}(\P(V_{\lambda}))=|\lambda|\,.
$
\end{itemize}

\subsubsection{The \emph{weight} truncations}\label{sssec-weight-filtration}
Let $p\in\NN$. Given an $\SSS_m$-module $W_m$, denote by $(W_m)_{\weight  > p}$ the sum of the irreducible factors of $W_m$ of weight $> p$.
\begin{itemize}
\item By Pieri's rule \ref{lemm-pieri}-(\ref{lemm-pieri-a}), if $\W:=\set \phi_m:W_m\to  W_{m+1}/_{m\in\NN}$ is an $\FI$-module, one has
$\phi_m\big((W_m)_{\weight >p}\big)\dans
\big((W_{m+1})_{\weight >p}\big)\,,$
in which case, the family 
$$\W_{\weight >p}:=\bigset\phi_m
(W_m)_{\weight >p}\to(W_m)_{\weight >p}/_{m\in\NN}$$
is a sub-$\FI$-module of $\W$. 

\item Let $\W_{\weight\leq p }:=\W/\W_{\weight> p }$. The short exact sequence
$$\0\to
\W_{\weight> p }
\to\W\to
\W_{\weight\leq p }
\to\0$$
is natural with respect to $\W$.

\end{itemize}

\begin{rema}The following are easy consequences of the  definitions.
\def\varlistskips{\itemsep0pt}\begin{enumerate}
\item The full subcategory  $\FImod_{\weight >p}$ of $\FI$-modules $\W$ such that the inclusion $\W_{\weight >p}\into\W$ is an isomorphism is an abelian subcategory.

\item  The full subcategory  $\FImod_{\weight \leq p}$ of $\FI$-modules $\W$ such that the quotient $\W\onto \W_{\weight \leq p}$ is an isomorphism is an abelian subcategory.

\item $\Ext^{i}_{\FI}\big(\FImod_{\weight> p},\FImod_{\weight\leq p}\big)=0$, for all $i\in\NN$.

\end{enumerate}
\end{rema}

\subsection{The $\FI$-module $\V_\lambda$}\label{Vlambda}Given a partition $\lambda$, consider the projective $\FI$-module $\P(V_{\lambda})$ of~\ref{ssec-proj}.  Lemma \ref{lemm-pieri}-(\ref{lemm-pieri-a}) says that for all $m\geq|\lambda|$ the smallest weight of the irreducible components of $\P(V_{\lambda})_{m}$ is exactly $|\ulambda|$ $({=}\weight(\lambda))$, which is the weight of a unique factor\,: the simple $\SSS_m$-module $V_{\ulambda[m]}$ with multiplicity $1$. The following proposition results from this simple observation and the weight filtration \ref{sssec-weight-filtration}.

\begin{prop}\label{prop-Vlambda}Let $\lambda$ be a nonempty partition.
\begin{enumerate}
\item 
\label{prop-Vlambda-a}
The terms of the quotient $\FI$-module
$$\V_{\lambda}:=\P(V_{\lambda})_{\weight\leq\weight(\lambda)}
=\bigset\phi_m:\V_{\lambda,m}\to \V_{\lambda,{m+1}}/_{m\in\NN}
\postskip1ex$$

\noindent are\hfil
$\V_{\lambda,m}=\begin{cases}\noalign{\kern-2pt}
\hfill\0\hfill\,,&\text{for all $m<|\lambda|$,}\\
V_{\ulambda[m]}\,,&\text{otherwise.}\\\noalign{\kern-2pt}
\end{cases}
\postskip1.75ex$

\medskip\noindent
The interior maps $\phi_m$ are injective, and are exhaustive for $m\geq |\lambda|$.

\item 
\label{prop-Vlambda-b}
If $\V'_{\lambda}:=\set\phi'_m:\V_{\lambda,m}\to \V_{\lambda,m+1}/_{m\in\NN}$ is an $\FI$-module such that $\phi'_{m}$ is injective and is exhaustive for all $m\geq N$, then $(\V'_\lambda)_{\geq N}$ and $(\V_\lambda)_{\geq N}$ (\cf\ref{comm-trunc}) are isomorphic $\FI$-modules for $m\geq N$.
\end{enumerate}
\end{prop}

\subsection{Representation stability of $\FI$ and of $\FB$ modules}\label{ssec-RS-modules}
\vskip-0.25ex\vskip0pt
\subsubsection{Representation stable $\FI$-modules}\label{sssec-RS-FI-modules}An $\FI$-module 
$\W=\set \phi_m:W_m\to W_{m+1}/_{m\in\NN}$ is said to be \emph{representation stable for $m\geq N$}, in short $\RS_{m\geq N}$, if the following conditions are satisfied.
\def\varlistskips{\def\theenumi{RS-\arabic{enumi}}\advance\itemindent0em
\advance\leftmargin2em
\itemsep2pt\parsep2pt\topsep2pt}\begin{enumerate}
\item
\label{ssec-RS-modules-RS-1}
The interior maps $\phi_{m}$ are injective, and are exhaustive for $m\geq N$.
\item
\label{ssec-RS-modules-RS-2}
For all $m\geq N$, we have a stable decomposition in simple modules
$$W_m\sim\bigoplusnl_{|\lambda|\leq N}V_{\ulambda[m]}{}^{\ng_{\ulambda}}\,,$$
where the  $\ng_{\lambda}$ are independent of $m\geq N$.
\end{enumerate}
We denote by $\rankFIRS(\W)$, the smallest such $N$, and we call it the \emph{rank of representation stability of $\W$.}

\noendpoint\begin{exas}\label{exas-weight}
\def\varlistskips{\itemsep1pt\parsep0pt\topsep1pt}\begin{enumerate}
\def\entree#1{\leavevmode\hbox to3.8cm{\,#1\hss}}
\item \entree{$\P(V_{\lambda})$ is $\RS_{m\geq|\lambda|+\lambda_1}$,} (lemma \ref{lemm-pieri}-(\ref{lemm-pieri-b})).
\item \label{exas-weight-b}
\entree{$\P(W_n)$ $\RS_{m\geq 2n}$,} 
 (consequence of lemma \ref{lemm-pieri}-(\ref{lemm-pieri-b})).
\item \entree{$\V_{\lambda}$ is $\RS_{m\geq |\lambda|}$,} (by definition). 
\end{enumerate}
\end{exas}

\subsubsection{Representation stable $\FB$-modules}\label{sssec-RS-FB-modules}An $\FB$-module 
$\W=\set W_m/_{m\in\NN}$ is said to be \emph{representation stable for $m\geq N$}, in short $\RS_{m\geq N}$, if the previous condition (\ref{ssec-RS-modules-RS-2}) is satisfied.

\smallskip
In that case, we say that $\W$ and $\bigoplus_{\lambda\in P}\V_{\lambda}^{\ng_\lambda}$ are \emph{asymptotically isomorphic as $\FB$-modules}, and we may write
$$\relax{\W_{\geq N}\sim\bigoplusnl_{|\lambda|\leq N}(\V_{\lambda})_{\geq N}^{\ng_\lambda}}\eqnnb\label{RS}$$

We denote by $\rankRS(\W)$, the smallest such $N$, and we call it the \emph{rank of representation stability of $\W$.} 

\noindent {\bf Warning. }For an $\FI$-module the stability rank \emph{as $\FB$-module} is, a priori, smaller than the stability rank as $\FI$-module. It is very easy to construct examples where it is strictly smaller, for example take any $\RS$ $\FI$-module and replace all  its interior maps by the zero map.

\noendpoint\begin{prop}\label{prop-RS-FB-modules}
\def\varlistskips{\itemsep1pt\parsep1pt\topsep0pt}\begin{enumerate}
\item\label{prop-RS-FB-modules-a}\hfil 
$\rankRS(\P(V_{\lambda}))=\rankFIRS(\P(V_{\lambda}))=|\lambda|+\lambda_1$.\vrule height12pt depth8pt width0pt
\item\label{prop-RS-FB-modules-b} If the trivial representation $\trivial_n$ is contained in $W_n$, then
$$\rankRS(\P(W_n))=\rankFIRS(\P(W_n))=2n\,.$$
\item\label{prop-RS-FB-modules-c}\hfil
$\rankRS(\V_{\lambda})=\rankFIRS(\V_{\lambda})=|\lambda|$. 
\vrule height12pt depth8pt width0pt
\item\label{prop-RS-FB-modules-d}
If $\W$ is an $\FB$-module such that $\Z\dans\W$, where
$\Z$ is any of the $\FI$-modules in the previous claims, then 
$$\preskip0.5ex
\rankRS(\W)\geq\rankRS(\Z)\,.$$
\end{enumerate}
\end{prop}
\demo\parskip2pt\mou
In (\ref{prop-RS-FB-modules-c}) there is nothing to be proved.
Let us denote by $\Z$ any of the other two $\FI$-modules in the claims.
We already gave the values of $\rankFIRS(\Z)$ in \ref{exas-weight}. Next, because $\rankRS(\Z)\leq \rankFIRS(\Z)$, the equality will follow if we show that for $m=\rankFIRS(\Z)$ the $\SSS_m$-module $Z_m$ contains an irreducible factor corresponding to a Young diagram having its two first lines of equal lengths\,:
$$\dimcell1.5mm
\defcellmacrowithin{1}{14}{\greyYNGtrue}
\putat{21mm}{1.4mm}{$\hyoung{1,1}$}
\putat{15.75mm}{0.8mm}{$\dots$}
\putat{15.75mm}{-0.5mm}{$\dots$}
V_{\hyoung{7,7,5,3}}
\hskip1cm\dans Z_m
\eqno(\diamond)$$
Indeed, if $\rankRS(\Z)<m$, we would have
$$Z_m\sim\bigoplusnl_{|\lambda|<m}V_{\ulambda[m]}{}^{\ng_{\ulambda}}\,,\eqno(\ast)$$
and the length of the first line of the Young diagram $\ulambda[m]$ in $(\ast)$, would be strictly greater than that of the others, \idest
$\ulambda[m]_{1}>\ulambda[m]_{2}$, which is contrary to our assumption (\cf\ref{lemm-pieri}-(\ref{lemm-pieri-b})).

Let us check $(\diamond)$ in our cases. In (\ref{prop-RS-FB-modules-a}) this results by Pieri's rule because $m=|\lambda|+\lambda_1$ is the smallest $m$ such that $\lambda$ is the socle of a diagram $\mu\vdash m$, which means that $\mu_1=\mu_2$. In (\ref{prop-RS-FB-modules-b}) the Young diagram corresponding to the trivial module $\trivial_{n}$ has only one line with $n$ boxes, and $\P(\trivial_n)_{2n}$ verifies~$(\diamond)$ for $m=2n$.

\noindent (\ref{prop-RS-FB-modules-c}) $W_m$ contains some irreducible $\SSS_m$-module of the form $(\diamond)$.
\enddemo

\section{Character polynomiality of $\FI$-modules}\label{sec-poly}

\vskip-1ex\vskip0pt\subsection{Character polynomiality}\label{ssec-poly}%
Let $\QQ_{\class}(\SSS_m)$ denote the algebra of \emph{rational class functions} on~$\SSS_m$, \idest functions $f:\SSS_m\to \QQ$ which are constant on each conjugacy class of~$\SSS_m$.
Denote by $\QQ[\cl X]$  the ring of polynomials with rational coefficients and countably many variables $X_1,X_2,\ldots,\ $ endowed with the  grading `$\degw$' that stipulates that\,:
$$\degw (X_i):=i\,.$$

\subsubsectionline{The \emph{weight\/} of a polynomial.}To avoid confusion in the sequel with the usual degree $\deg(P)$ of a polynomial $P\in \QQ[\cl X]$ which stipulates that $\deg (X_i)=\nobreak1$, we will call $\degw(P)$ the \emph{weight} of $P$.

%

\begin{prop}\label{prop-poly}
Denote by $X_{m,i}:\SSS_m\to\NN$ the class function which assigns to $g\in \SSS_m$, the number $X_{m,i}(g)$ of $i$-cycles in the decomposition of $g$ as product of disjoint cycles  in $\SSS_m$.
\begin{enumerate}
\item
\label{prop-poly-a}
The map
$$
\preskip-10pt\mathalign{\rho_m:&\QQ[\cl X]&\too& \QQ_{\class}(\SSS_m)\hfill\\
&\hfill X_i&\mapstoo&(g\mapsto X_{m,i}(g))}\eqnnb\label{rho-m}$$
is an homomorphism of algebras whose kernel contains the polynomials
$$\mathrigid2mu
\big(X_1+2X_2+\cdots+mX_m-m\big)
\text{\quad and \quad}
\big(X_i(X_i-1)\cdots\big(X_i-\big\lfloor m/i\big\rfloor\big)\big)\,.
\eqnnb\label{rho-m-ker}$$
The restriction 
$$\rho_m:\QQ[X_1,\ldots,X_{m-1}]\onto\QQ_{\class}(\SSS_m)$$
is surjective. In particular, the characters of $\SSS_m$ are represented by polynomials  with rational coefficients and variables $X_1,\ldots, X_{m-1}$. 
\item
\label{prop-poly-b}
For $n\leq m$ and $g\in\SSS_n$, we have
$$
\preskip-0.8ex
\postskip-1ex
\mathalign{
({\rm i})\ &\rho_{m}(X_1)(\iota g)&=&\rho_{n}(X_1)(g)+(m-n)\,,\\\noalign{\kern2pt}
({\rm ii})\ &\rho_{m}(X_i)(\iota g)&=&\rho_{n}(X_i)(g)\,,\ \forall i>1\,.\hfill
}
\quad
\varxymatrixc{3mm}{@R=1mm@C=1.2cm}{
\vrule depth4pt width0pt \SSS_n\ar[rd]|{\ \rho_n\ }\ardans[dd]_{\iota}\\
&A\\
\SSS_m\ar[ru]|{\ \rho_m\ }\\
}
$$\vskip0mm
\end{enumerate}\end{prop}
\def\demoname{\it Hint. }
\demo (\ref{prop-poly-a}) That the polynomials \ref{rho-m-ker} belong to $\ker(\rho_m)$ is clear. Next, to see that $\rho_m$ is surjective,
we need only show that the characteristic function of a conjugacy class of $\SSS_m$ can be realized as a polynomial in $X_1,\ldots, X_m$. 

\scaledisplayskips90/100

\noindent For $k\in\iii[1,m]$, let
$R_{k}(Z):=Z(Z-1)\cdots\widehat{(Z-k)}\cdots\big(Z-m)
$
and consider
$$D_{k}(Z):=R_{k}(Z)/R_{k}(k)\in\QQ[Z]\,.$$
This polynomial has the property that
$$
\postskip4pt
\rho_m(D_{k}(X_i)(g)=\left\{\mathalign{1\,,&\text{ if $X_i(g)=k$\,,}\\
0\,,&\text{ otherwise.}\hfill 
}\right.$$
So that, if $\sum_i i\,\ng_i=m$, we get
$$
\preskip4pt
\rho_m(D_{\ng_1}(X_1) D_{\ng_2}(X_2)\cdots D_{\ng_m}(X_m))(g)
=
\left\{\mathalign{1\,,&\text{ if $g$ is of type $(1^{\ng_1},\ldots, m^{\ng_m})$}\\
0\,,&\text{ otherwise.}\hfill 
}\right.
$$
(\ref{prop-poly-b}) is clear.
\enddemo

\killline
\def\vartext{Convention}
\begin{var}In order to alleviate notations, 
 we will write $X_i(g)$ for $X_{m,i}(g)$  despite the possible ambiguity in `$X_1(g)$' \emph{(\cf \ref{prop-poly}-(\ref{prop-poly-b}-i))}.
\end{var}

\begin{defi}\label{defi-polynomial-char}An $\FB$-module $\W:=\set W_m/_{m\in\NN}$ is said to have a \emph{polynomial character for $m\geq N$}, in short $\PC_{m\geq N}$, if there exists a polynomial $P_{\W}\in \QQ[\cl X]$ such that
$$\chi_{W_m}=\rho_{m}(P_{\W})\,,\mrlap{\quad\forall m\geq N\,.}$$
We denote by $\rankPC(\W)$, the smallest such $N$, and we call it the \emph{rank of character polynomiality of $\W$.}
\end{defi}

\begin{prop}[and definition]\label{prop-polynomial-char}The polynomial $P_\W$ that asymptotically represents the characters of the terms $W_m$ of  an $\FB$-module $\W:=\set W_m/_{m\in\NN}$ is
unique, it is called
the \emph{polynomial character of the $\FB$-module $\W$}.
\end{prop}
\demo If $P'_{\W}$ were another polynomial representing $\chi_{W_m}$ for $m\ggg 0$, the difference $Q:=P_{\W}-P'_{\W}$, that we may assume to belong to $\QQ[X_1,\ldots,X_N]$, would be a polynomial representing the zero class function for all $m\ggg N$. 

If $Q$ is not the null polynomial, we  can write it as a polynomial in $X_N$ with coefficients $Q_i$ in $\QQ[X_1,\ldots,X_{N-1}]$\,:
$$Q=Q_0 + Q_1 X_N+ Q_2 X_N^2+\cdots+ Q_r X_N^r\,,\text{ and $Q_r\not=0$.}\eqno(\ast)$$
Now, for any family $\cl a:=\set a_1,\ldots,a_{N-1}/\dans\NN$, and any
$i\in\NN$, it is easy to find $m_i\ggg N$ and $g_i\in \SSS_{m_i}$ such that $X_1(g_i)=a_1$, \dots, $X_{N-1}(g_i)=a_{N-1}$ and  $X_{N}(g_i)>i$. In that case, the polynomial in one variable $Q(\cl a,X_N)$ has infinitely many roots and is, therefore, the null polynomial. In particular, in  the decomposition $(\ast)$, we would have \hbox{$Q_r(\cl a)=0$} for all choices of $\cl a$, but this implies that $Q_r$ is the null polynomial in $\QQ[X_1,\ldots,X_{N-1}]$, contrary to its definition. The polynomial $Q$ must therefore be the null polynomial.
\enddemo

\subsection{Frobenius character formula}\label{ssec-Frobenius-char}Given a partition $\lambda:=\parti(\lambda_1,\dots,\lambda_{\ell})\vdash m$, the celebrated Frobenius formula computes the character $\chi_{V_\lambda}$ of the simple $\SSS_m$-module $V_{\lambda}$. The important point for us about this  formula is that it gives an expression of $\chi_{V_\lambda}$ as a polynomial \emph{depending only on the socle $\ulambda$} of $\lambda$. Consequently, the same polynomial expresses the characters of all the terms of the $\FI$-module $\V_{\lambda}=\set V_{\lambda,m}/_{m}$, for $m\geq |\lambda|$, which self-explains the property of character polynomiality of $\V_{\lambda}$.

\subsubsection{Frobenius polynomial for $\chi_{V_{\lambda}}$}\label{sssec-Frobenius-char}Following Macdonald in his book \cite{mac-livre} (ex. I.7.14, p. 122), let $y:=\set y_1,\ldots,y_\ell/$ be a set of $\ell$  variables, where $\ell:=\ell(\lambda)$. The \emph{discriminant of $y$} is the antisymmetric homogeneous polynomial
$$\Delta(y):=\prodnl_{i<j}(y_i-y_j)\,.$$
For  $d\in\NN$, the \emph{$d$-power sum of $y$} is the symmetric homogeneous polynomial
$$P_{d}(y):=y_1^{d}+\cdots+y_\ell^{d}\,.$$
The value $\chi_{V_{\lambda}}(g)$ for $g\in \SSS_{m}$, is the coefficient of the monomial 
$$
y_1^{\lambda_1+(\ell-1)}
y_2^{\lambda_2+(\ell-2)}
y_3^{\lambda_3+(\ell-3)}\cdots
y_\ell^{\lambda_\ell}\,,
$$
in the development of the product
$$\Delta(y)\,\Big(\prodnl_{d\geq1}P_{d}(y)^{X_d(g)}\Big)\,.\eqnnb\label{eqn-frobenius-formula}$$
This coefficient, denoted by $\Xg_{\lambda}$, is a polynomial in $\QQ[\cl X]$, we call it \emph{the Frobenius polynomial for $\chi_{V_{\lambda}}$}.

\begin{prop}\label{prop-Frobenius-char}The Frobenius polynomial $\Xg_{\lambda}$ for $\chi_{V_{\lambda}}$ depends only on the socle $\ulambda$ of $\lambda$, and belongs to the ring
$\QQ[X_1,\ldots,X_{\lambda_2+\ell-2}]$. Its weight is\,:
$$\degw(\Xg_\lambda)=\weight(\lambda)\,.$$
The characters of\/ $\SSS_m$ can be represented by polynomials in $\QQ[X_1,\ldots, X_{m-1}]$ of weights $\leq m-1$.
\end{prop}
\demo
Because the polynomial \ref{eqn-frobenius-formula} is homogeneous, we can make $y_1=1$ without loosing information. In that case, $\chi_{V_\lambda}(g)$ is the coefficient in 
$$
y_2^{\lambda_2+(\ell-2)}
y_3^{\lambda_3+(\ell-3)}\cdots
y_\ell^{\lambda_\ell}\,,
$$
after the development of the product
$$\Delta(\tilde y)\,\Big(\prodnl_{j>1}(1-y_j)\Big)\,
\Big(\prodnl_{d\geq1}(1+P_{d}(\tilde y))^{X_d(g)}\Big)\eqno(\ddagger\ddagger)$$
where $\tilde y:=\set y_2,\ldots,y_\ell/$.
But, in this product the first factor $\Delta(\tilde y)$ is already homogeneous of total degree $(\ell-2)+(\ell-3)+\cdots$, so that we must seek, in the development of the remaining factors, terms whose total degree is bounded by $|\ulambda|=\lambda_2+\cdots+\lambda_\ell$. But then, since we have
$$\big(1+P_d\big)^{X_d}=1+\binom{X_d}1P_d+\binom{X_d}2P_d^{2}
+\binom{X_d}3P_d^{3}+\cdots$$
and because $\totdeg(P_d^{a})=ad$, we conclude that
\begin{itemize}
\item If $d>\lambda_2+\ell-2$, the factor $\big(1+P_d\big)^{X_d}$  contributes only to $\chi_{V_\lambda}$ with its term $1^{X_d}$, so can be ignored. The product symbol $\prod_{d\geq 1}$ in $(\ddagger\ddagger)$ can therefore be replaced by $\prod_{d=1}^{\lambda_2+\ell-2}$.

\item The coefficient $\binom {X_d}j$ is a polynomial of degree $j$ in $X_d$ and appears attached to monomials in $\tilde y$ of total degree $jd$. We can thus infer that, after development, the expression of $\chi_{V_{\lambda}}$ is a polynomial in $X_1,\ldots,X_{\lambda_2+\ell-2}$ of \emph{weight\/} $|\ulambda|=\weight(\lambda)$.\enddemobox
\end{itemize}
\enddemo

\noindent The following corollary of proposition \ref{prop-Frobenius-char} is now immediate from the definition of representation stable $\FB$-modules \ref{sssec-RS-FB-modules}.\killline

\noendpoint\begin{coro}\label{coro-Frobenius-char}If $\W$ is an $\FB$-module, then 
$$
\rankPC(\W)\leq\rankRS(\W)\,.
$$
\end{coro}
\demo Left to the reader.
\enddemo

\subsection{Basic examples of character polynomiality of $\FI$-modules}\label{ssec-basic-poly}

\subsubsection{The $\ell$-cycles of $\iii[1,m]$}\label{sssec-cycles}

Given $m,\ell\in\NN$, we denote by $\iGamma^m_\ell$ the set of
$\ell$-tuples $(i_1,\ldots,i_\ell)$ of pairwise distinct elements of $\iii[1,m]$ modulo cyclic permutation, \idest such that
$$\def\sim{=}(i_1,\ldots,i_\ell)\sim(i_2,\ldots,i_\ell,i_1)\sim(i_3,\ldots,i_\ell,i_1,i_2)\sim\cdots$$
The symmetric group $\SSS_m$ acts on $\varGamma^m_\ell$ by
$$g\cdot(i_1,\ldots,i_\ell)=(g(i_1),g(i_2),\ldots,g(i_\ell))\,.\eqnnb\label{Sm->iGamma}$$
The elements of $\iGamma^m_\ell$ are called \emph{the $\ell$-cycles of $\iii[1,m]$}.

\noendpoint\begin{comms}\label{comms-cycles}
\def\varlistskips{\listparindent0pt\itemsep2pt\mou\parsep0pt
}\begin{enumerate}
\item 
\label{comms-cycles-a}Given $g\in\SSS_m$, the set $\iii[1,m]$ is decomposed in $\langle g\rangle$-orbits, each of which can be endowed with a cyclic order defined by $g$. For example, if $x\in\iii[1,m]$, we may consider $(x\to g(x)\to g^2(x)\to\cdots\to x)$, which gives the well-known \emph{decomposition of $g\in\SSS_m$ as product of disjoint cycles.}

\item 
\label{comms-cycles-b}For $m\geq\ell$, there is a difference between the cases $\ell=1$ and $\ell>1$. 

\smallskip
For all $\ell\geq1$, define the \emph{support} of an $\ell$-cycle\/ $\gamma:=(i_1,\ldots,i_\ell)$ the set $\supp|\gamma|$  of its coordinates\,: 
$$\supp|\gamma|:=\set i_1,\ldots,i_\ell/\dans\iii[1,m]\,.$$ 

Then, define $\itilde\gamma\in\SSS_m$ by\,:
$$\itilde\gamma:=\mycases{
\itilde \gamma(i_j)&:=&i_{j+1\mod\ell}\,,&\text{ for $i_j\in\supp|\gamma|$,}\\
\itilde\gamma(x)&:=&x\,,\hfill&\text{ for $x\not\in\supp|\gamma|$.}}$$ 

\def\varlistskips{\itemsep2pt\mou\parsep0pt\leftmargin\parindent
\listparindent0pt
}\begin{enumerate}
\item
\label{comms-cycles-bi} For $\ell=1$, we have $\iGamma_1^m=\iii[1,m]$ 
and the action \ref{Sm->iGamma} of $\SSS_m$ on~$\iGamma_1^m$
is the standard action of $\SSS_m$ on $\iii[1,m]$. The map ${(\widetilde{\_})}:\iGamma_{1}^m\to\SSS_m$ is uninteresting, as it is the constant map $\gamma\mapsto\id_{\iii[1,m]}$.

\item
\label{comms-cycles-bii}
For $\ell>1$, the map
${(\widetilde{\_})}:\iGamma_\ell^m\dans\SSS_m$ is \emph{injective}, and the action \ref{Sm->iGamma} of $\SSS_m$ on~$\iGamma_\ell^m$ is the conjugation action of $\SSS^m$ on itself, \idest\unskip\,:
$$\halfdisplayskips
\widetilde{g\cdot\gamma}=g\,\itilde\gamma\, g^{-1}\,.$$

\smallskip
We will identify $\gamma$ and $\itilde\gamma$ if no confusion is likely to arise. In this sense, for $g\in\SSS_m$, the set of fixed points
 $(\iGamma_{\ell}^{m})^{g}:=\bigset \gamma\in\iGamma_{\ell}^{m}\mid g\cdot\gamma=\gamma/$ is simply\,:
$$(\iGamma_{\ell}^{m})^{g}=\bigset \text{$\ell$-cycles\ }\gamma\in\SSS_m\mid g\gamma=\gamma g/\,.\eqnnb\label{g-fixed-cycles}$$

\item 
\label{comms-cycles-biii}For all $m,\ell>0$, we have\,:
$$|\iGamma_\ell^m|=\frac{m(m-1)\cdots(m-(\ell{-}1))}\ell\,\cdot$$
\end{enumerate}\end{enumerate}
\end{comms}

\subsubsection{The $\FI$-modules $\IE_{\nu}$}
For $m,\ell>0$, let $\IE^{m}_{\ell}$ be the vector space spanned by the set $\iGamma^m_\ell$, \idest 
$$\IE^{m}_{\ell}:=\bigoplus\nolimits_{\gamma\in \iGamma^m_\ell} \QQ\cdot \gamma\,.$$
Endow it with the linear action of $\SSS_m$ induced by its action on the basis $\iGamma_\ell^m$.

\noindent Notice that\,:
\begin{itemize}
\mynobreak\nobreak
\item
\leavevmode \vrule height13pt width0pt  $\IE_\ell^m=\0$, for all $m<\ell$. 

\item For all $m\geq n$, the set
$\iGamma_{\ell}^{n}$ is a subset of $\iGamma_{\ell}^{m}$, which is invariant under the action of $\1_n\sqtimes\SSS_{m-n}$, so that the natural inclusions $\iGamma_{\ell}^{m}\dans \iGamma_{\ell}^{m+1}$ 
 induce the interior maps (clearly injective) of an $\FI$-module 
$$\displayboxit{\IE_\ell:=\bigset \phi(\IE_\ell)_m:\IE_\ell^m\into \IE_\ell^{m+1}/}$$
where $\phi(\IE_\ell)_m$ is exhaustive for $m\geq\ell$, since the natural map between orbit spaces: $\SSS_m\backslash\iGamma_{\ell}^{m}\to 
\SSS_{m+1}\backslash\iGamma_{\ell}^{m+1}$ is bijective for $m\geq\ell$.
\end{itemize}
\medskip\noindent
More generally, if $\nu:=(1^{\ng_1},2^{\ng_2},\ldots,N^{\ng_N})=(\nu_1\geq\nu_2\geq\cdots\geq\nu_\ell)
$ is a nonempty partition, let
$$\let\EE\IE
\mathalign{\EE_\nu^{m}&:=&
(\EE_{1}^{m})^{\otimes \ng_1}\otimes
(\EE_{2}^{m})^{\otimes \ng_2}\otimes
\cdots\otimes
(\EE_{N}^{m})^{\otimes \ng_N}\,,\hfill\\
&:=&\IE_{\nu_1}\otimes\IE_{\nu_2}\otimes\cdots\otimes
\IE_{\nu_\ell}\,,\hfill}$$
and define $\phi(\IE_{\nu})_m:\IE_\nu^m\to \IE_\nu^{m+1}$ as
the tensor product of interior maps 
$$
\mathalign{\phi(\IE_{\nu})_m&:=&\phi(\IE_1)^{\otimes\ng_1}
\otimes\cdots\otimes
\phi(\IE_N)^{\otimes\ng_N}\,,\\
&:=&\phi(\IE_{\nu_1})
\otimes\cdots\otimes
\phi(\IE_{\nu_\ell})\,.\hfill}
$$
The family 
$$\postskip1.em
\displayboxit{\IE_\nu:=\bigset \phi(\IE_{\nu})_m:\IE_\nu^m\into \IE_\nu^{m+1}/_m}$$
is an $\FI$-module with all interior maps injective, and exhaustive for $m\geq \nu_1$.

\noendpoint\begin{prop}\label{prop-basic-CP}
\begin{enumerate}
\item
\label{prop-basic-CP-a}
For $m,\ell\in\NN$, the character $\chi_{\IE^{m}_{\ell}}$ is expressed by the following polynomial of $\QQ[X_1,\ldots, X_\ell]$ of weight $\ell$ and independent of $m\in\NN$\,:
$$\def\le{{\ell/ e}}
\def\le{{d}}
\Eg_\ell=\phi(\ell)\,X_\ell+
\sumnl_{ed=\ell\,,\; e\not=1\;}
\phi(d)\,\frac{d^{e}}{\ell}\,
X_{\le}(X_{\le}-1)\cdots\big(X_{\le}-(e{-}1)\big)\,,\eqno(\ast_{\ell})$$
where $\phi$ is the Euler's totient function.

\smallskip
\noindent The $\FI$-module 
$\IE_{\ell}:=\set \IE_{\ell}^{m}\into\IE_{\ell}^{m+1}/_{m}$
has the following ranks
$$\rankPC(\IE_{\ell})=0\text{\quad and\quad }
\rankRS(\IE_{\ell})=\rankFIRS(\IE_{\ell})=2\ell\,.$$
\item\label{prop-basic-CP-b}Given an ordered sequence $\cl Z:=(Z_1,\ldots,Z_N)$ of polynomials of\/ $\QQ[\cl X]$, and $d\in\nobreak\NN$, denote by\/ $\QQ^{\leq d}[Z_1,\ldots,Z_N]$ the subspace of polynomials of weight~$\leq \nobreak d$ relative to $\cl Z$, \idest the subspace spanned by the elements $\def\ng{a}Z_1^{\ng_1}\cdots Z_N^{\ng_N}$ where $
\def\ng{a}\sum_i i\,\ng_i\leq d$. Then, for all $N\in\NN$, 
the inclusion\,:
$$
\QQ^{\leq d}[\Eg_1,\ldots,\Eg_N]\dans \QQ^{\leq d}[X_1,\ldots,X_N]$$
 is an equality.

\item\label{prop-basic-CP-c}For every nonempty partition $\nu:=(1^{\ng_1},2^{\ng_2},\ldots,N^{\ng_N})$,
we have
\begin{enumerate}
\item 
\label{prop-basic-CP-ci}
$\weight(\IE_{\nu}^{m})\leq\min\set m,|\nu|/$.
\item
\label{prop-basic-CP-cii} The $\FI$-module $\IE_\nu:=\set\IE_\nu^m\into\IE_\nu^{m+1}/_m$ 
has the following ranks
$$\rankPC(\IE_{\nu})=0\text{\quad and\quad }
\rankRS(\IE_{\nu})=\rankFIRS(\IE_{\nu})=2|\nu|\,.$$

\end{enumerate}
\end{enumerate}
\end{prop}
\demo (\ref{prop-basic-CP-a}) Since the linear action of $g\in\SSS_m$ is induced by its action on the basis $\iGamma^m_\ell\dans\IE_{\ell}^{m}$, the trace $\chi_{\IE_\ell^{m}}(g)$ is the cardinality of  $(\iGamma^m_\ell)^{g}$ (\ref{comms-cycles}-(\ref{comms-cycles-bii})). 

When $m<\ell$, the set $\iGamma^{m}_{\ell}$ is empty and $\IE_\ell^{m}=0$.
And, since in $\SSS_m$ there is no permutation $g$ such that $\rho_{m}(X_d)(g)\geq \ell/d$, we have  $\rho_{m}(\Eg_{\ell})=0$, which
states (\ref{prop-basic-CP-a}) when $m<\ell$. 

When $m\geq\ell$, following \ref{comms-cycles}-(\ref{comms-cycles-b}), we  consider two cases\,:

-- {\boldmath$\ell=1$. }Then $\iGamma_1^m=\iii[1,m]$, $\chi_{\IE_1^m}(g)=|\iii[1,m]^{g}|=X_1(g)$, and (\ref{prop-basic-CP-a}) is obvious.

-- {\boldmath$\ell>1$. }The set $(\iGamma^m_\ell)^{g}$ identifies, after \ref{comms-cycles}-(\ref{comms-cycles-bii})-\ref{g-fixed-cycles}, with the set of $\ell$-cycles $\gamma\in \SSS_m$ such that
$g\gamma=\gamma g$. As a consequence, $g\supp|\gamma|=\supp|\gamma|$ and the set $\supp|\gamma|$ appears endowed with two actions commuting to each other. But then, since the action of $\gamma$ on $\supp|\gamma|$ is transitive, the $\langle g\rangle$-orbits share the same cardinality. Denote it by $d$ and set $e:=\ell/d$. Then, 
$g^{\kappa}=\gamma^{e}$,  for some $\kappa\in\iii[1,d]$, relatively prime to $d$. 
These observations are gathered in the following image:
$$
\includegraphics{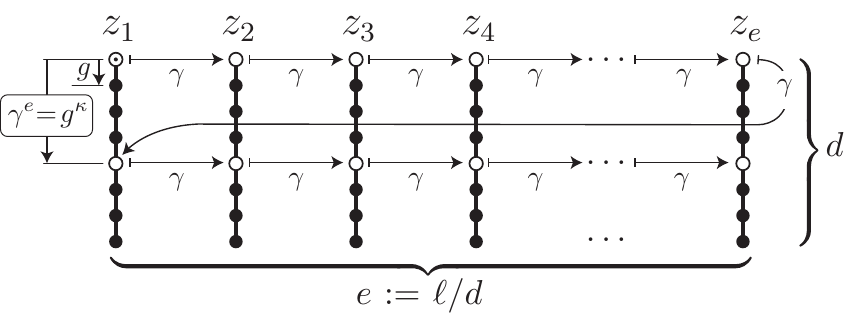}
\postskip4pt
$$
If we choose the point `\raise1pt\hbox{$\scriptstyle\odot$}' to stand for the smallest coordinate in~$\gamma$, this graphic representation becomes one-to-one,
and thus helpful for the computation of the cardinality of the set $(\iGamma^{m}_{\ell})^{g}$.
Indeed, it can be used to reconstruct $(\iGamma^{m}_{\ell})^{g}$ by following the next steps.
\def\varlistskips{\parsep0pt\topsep4pt\def\theenumi{\arabic{enumi}}}\begin{enumerate}
\item Take an ordered sequence $(z_1,\ldots,z_e)$ of $d$-cycles in $g$, among the
$$X_d(g)^{\usp e}:=X_d(g)\big(X_d(g)-1\big)\cdots
\big(X_d(g)-(e{-}1)\big)\,.
\postskip2pt
$$
possible such sequences.
\item Shift to the first position the $d$-cycle containing the smallest coordinate in $\iii[1,m]$. This means that $X_d(g)^{\usp e}$ has to be divided by $e$. 
\item Fix the vertical orders of $z_2,\ldots,z_e$. There are $d^{e-1}$ such orders. 
\item Fix the number $\kappa$. There are $\phi(d)$ possibilities.
\end{enumerate}
Therefore,
$$
\mathalign{\big|(\iGamma_{\ell}^{m})^g\big|
&=&
\sum_{de=\ell}\,\phi(d)\,\frac{d^{e}}\ell\,
X_d(g)^{\usp e}\,.\hfill
}
$$

\noindent In both cases, $\ell=1$ and $\ell>1$, the previous expression is independent of $m\in\NN$, which means that $\rankPC(\IE_\ell)=0$.

\smallskip
To estimate $\rankRS(\IE_{\ell})$, observe that, since the action of $\SSS_m$ on $\iGamma_{\ell}^m$ is transitive, $\iGamma^{m}_{\ell}=\SSS_m\cdot(1,\ldots,\ell)$, so that, denoting $\gamma:=(1,\ldots,\ell)$, we have an isomorphism of $\SSS_m$-spaces\,: 
$$\iGamma_\ell^m=\SSS_m\cdot\gamma\simeq\SSS_m\times_{\langle\gamma\rangle\sqtimes\SSS_{m-\ell}}\set\mathord{\rm pt}/\,,$$ 
because $\Stab_{\SSS_m}(\gamma)=\langle\gamma\rangle\sqtimes\SSS_{m-\ell}$.
Therefore, 
$$\IE_\ell^m=\ind_{\SSS_\ell\sqtimes\SSS_{m-\ell}}^{\SSS_m}
\big((\ind^{\SSS_\ell}_{\langle\gamma\rangle}\QQ)\sqtimes \trivial_{m-\ell}\big)\,.$$
and $\IE_\ell$ is the projective $\FI$-module $\P\big(\ind^{\SSS_\ell}_{\langle\gamma\rangle}\QQ\big)$ for which we have already shown that it is
$\RS_{m\geq2\ell}$ (\cf\ref{exas-weight}-(\ref{exas-weight-b})).

We can be a more precise observing that the simple $\SSS_\ell$-module~$\trivial_{\ell}$ appears with multiplicity $1$ in $\ind^{\SSS_\ell}_{\langle\gamma\rangle}\QQ$, which is true since\,:
$$
\Hom_{\SSS_\ell}(\ind^{\SSS_\ell}_{\langle\gamma\rangle}\QQ,\trivial_{\ell})=
\Hom_{\langle\gamma\rangle}(\QQ,\QQ)=\QQ\,.$$
Thus, $\IE_\ell$ contains
$\P(\trivial_\ell)=\P(V_{0[\ell]})$  as the sub-$\FI$-module, and, thanks to
\ref{prop-RS-FB-modules}-(\ref{prop-RS-FB-modules-d}),
$$\rankRS(\IE_\ell)=\rankFIRS(\IE_\ell)=2\ell\,.$$

\smallskip\noindent(\ref{prop-basic-CP-b}) 
In (\ref{prop-basic-CP-a}), 
the equation $(\ast_{\ell})$ shows that $\Eg_{\ell}=\phi(\ell)\,X_\ell$ modulo a polynomial of weight $\ell$ in the variables $X_d$ for $d|\ell$ and $d\not=\ell$. (In particular, $\Eg_{1}=X_1$.)
Now, given $N\geq1$, the system of equations $(\ast_{\ell})$, for $\ell\leq N$, can be inverted by recursively introducing some polynomials $Q_{\ell}^{m}(Z_1,\ldots,Z_\ell)\in\QQ[\cl Z]$ of weight~$\ell$, such that
$$X_{\ell}=Q^{m}_{\ell}(\Eg_1,\ldots,\Eg_\ell)\,,\relax{\quad\forall \ell\leq N\,.}\postskip4pt$$
The equality,
$$\QQ[X_1,\ldots,X_N]=\QQ[\Eg_1,\ldots,\Eg_N]\,,\postskip4pt$$
but also,
$$\QQ^{\leq d}[X_1,\ldots,X_N]=\QQ^{\leq d}[\Eg_1,\ldots,\Eg_N]\,,\quad\forall d\in\NN\,,$$
then follows easily.

\medskip\noindent(\ref{prop-basic-CP-c}) Since the action of $\SSS_m$ on $\IE_\nu^m:=\IE_{\nu_1}^m\otimes\cdots\otimes
\IE_{\nu_\ell}^m$ is induced by its component-wise action on the basis 
$
\iGamma^m_{\nu}:=\iGamma^m_{\nu_1}\times\iGamma^m_{\nu_2}\times\cdots\times\iGamma^m_{\nu_\ell}$, the structure of $\QQ[\SSS_m]$-module of $\IE_\nu^m$ follows from the structure of $\SSS_m$-space of $\iGamma^m_{\nu}$. Given $\cl \gamma:=(\gamma_1,\ldots,\gamma_\ell)\in\iGamma^m_{\nu}$, we have
$$\SSS_m\cdot\cl \gamma=
\SSS_m/\Stab_{\SSS_m}(\cl \gamma)=
\SSS_m\times_{\SSS_{\supp|\cl\gamma|}}\big(
\SSS_{\supp|\cl\gamma|}/\Stab_{\SSS_{\supp|\cl\gamma|}}(\cl \gamma)
\big)
\eqno(\ddagger)$$
where we set $\supp|\cl\gamma|:=\supp|\gamma_1|\cup\cdots\cup\supp|\gamma_\ell|$.

As a consequence the $\QQ[\SSS_m]$-module $\IE_\nu^{m}$ is a direct sum of induced modules of the form:
$$\ind^{\SSS_m}_{\SSS_{\supp|\cl\gamma|}\sqtimes\SSS_{\complement\supp|\cl\gamma|}}W_n\sqtimes \trivial_{m-n}\,,$$
where we can choose to have $\supp|\cl\gamma|:=\iii[1,n]$, for 
 $n:=|\supp|\cl\gamma||\leq\min\set m,|\nu|/$, in which case
 $$W_n:=\QQ[\SSS_{n}/\Stab_{\SSS_{n}}(\cl \gamma)]
 =\ind^{\SSS_{n}}_{\Stab_{\SSS_{n}}(\cl \gamma)}\trivial_{n}
 \,.\eqno(\ddagger\ddagger)
$$

When $m\geq|\nu|$, the corollary
\ref{lemm-pieri}-(\ref{lemm-pieri-a}) to Pieri's rule, gives the inequality
$$\weight(\nu)\leq\weight(\IE_{\nu}^m)\leq|\nu|\,,
$$
and (\ref{prop-basic-CP-ci}) follows.

\smallskip\noindent
(\ref{prop-basic-CP-cii}) Since $\chi_{\IE_\nu^m}=(\chi_{\IE_1^m})^{\ng_1}\cdots
(\chi_{\IE_N^m})^{\ng_N}$ and since $\rankPC(\IE_\ell)=0$, after (\ref{prop-basic-CP-a}), it follows immediately that $\rankPC(\IE_\nu)=0$. 

To finish,  we notice that, on the one hand, the number of terms of the form~$(\ddagger)$ is exactly de number of $\SSS_m$-orbits in $\iGamma_{\nu}^m$, and, on the other hand, the sequence of inclusions $\SSS_m\backslash\iGamma_{\nu}^m\dans \SSS_{m+1}\backslash\iGamma_{\nu}^{m+1}$ is \emph{strict increasing} until it stabilizes for $m\geq |\nu|$. 
As a consequence, the truncated $\FI$-module $(\IE_\nu)_{\geq |\nu|}$ is  isomorphic to the sum of truncated projective $\FI$-modules $\P(W_n)_{\geq|\nu|}$, all of which being $\RS_{m\geq 2|\nu|}$, after \ref{exas-weight}-(\ref{exas-weight-b}).
We therefore have proved that
$$\preskip6pt\rankRS(\IE_{\nu})\leq\rankFIRS(\IE_{\nu})\leq 2|\nu|\,.$$
To see that these are equalities, it suffices, by \ref{prop-RS-FB-modules}, to show that, for $m=|\nu|$, one has\,:
$$\P(\trivial_{m})\dans\IE_{\nu}\,.\postskip1.5ex\eqno(\dagger)$$
But, this is clear if we choose $\cl\gamma\in
\iGamma_{\nu}^{m}$ such that $|\supp|\cl\gamma||=m$ (as in $(\ddagger\ddagger)$), in which case $(\dagger)$ follows, proceeding as in (\ref{prop-basic-CP-a}), from the inclusion:
$$\trivial_m\dans\ind^{\SSS_{m}}_{\Stab_{\SSS_{m}}(\supp|\cl\gamma|)}(\trivial_{m})\,.$$
\vskip-2.3em
\enddemo

\section{Polynomial weight vs. Representation weight}
\vskip-0.5ex\vskip0pt\subsection{On the minimal polynomial weight of a character}The Frobenius polynomial 
$\Xg_\lambda\in\QQ[\cl X]$ is of weight~$\weight(\lambda)$, and there are infinitely many polynomials $P_{\lambda}\in\QQ[X_1,\ldots,X_{m-1}]$ such that $\rho_m(P_{\lambda})=\rho_m(\Xg_{\lambda})=\chi_{V_{\lambda}}$. In this section, we advance the study of the relationship between weights of polynomials and weights of representations. We will see that $\weight(W_m)$ is the lowest possible weight for a polynomial $P_{W_m}\in\QQ[X_1,\ldots,X_{m-1}]$ verifying $\rho_m(P_{W_m})=\chi_{W_m}$. 
For example, the polynomial $\Xg_{W_m}:=\sum_\lambda\ng_\lambda\Xg_\lambda$, where $\Xg_{\lambda}$ is the Frobenius polynomial of an irreducible factor $V_\lambda$  of multiplicity~$\ng_\lambda$ in~$W_m$. Moreover, we will see that this polynomial $\Xg_{W_m}$ is the only one satisfying these weights conditions,
 if and only if $\weight(W_m)\leq m/2$.

\noendpoint\begin{theo}
\label{theo-degree-poly}
\def\varlistskips{\itemsep4pt\parsep0pt\topsep0pt}\begin{enumerate}
\item\label{theo-degree-poly-a}
For $d,m\in\NN$, denote by $\rho_{m|d}$ the restriction 
of $\rho_{m}:\QQ[\cl X]\to\QQ_{\class}(\SSS_m)$
to the subspace $\QQ^{\leq d}[\cl X]$ of polynomials of weight $\leq d$.
The image of the map
$$\postskip0pt
\rho_{m|d}:\big(\QQ^{\leq d}[\cl X]\big)\to
\QQ_{\class}(\SSS_m)
$$
is the subspace
$$
\QQ_{\class}^{\leq d}(\SSS_m)=
\big\langle
\,\chi_{V_{\lambda}}\mathbin{\big|}\big(|\lambda|=m\big)\mathbin \&\big(\weight(\lambda)\leq d\big)
\big\rangle_{\QQ}\,,
\eqnnb\label{Shur-ortho}$$
and $\ker(\rho_{m|d})=0$ if and only if $d\leq m/2$.
\item
\label{theo-degree-poly-b}
Let $P\in \QQ[\cl X]$.
\def\varlistskips{\itemsep0pt\parsep2pt\topsep0pt}\begin{enumerate}
\item
\label{theo-degree-poly-bi} The scalar product
$\scalar1{P}_{\SSS_m}$
is constant for $m\geq\degw (P)$.
\item
\label{theo-degree-poly-bii}
The scalar product
$\scalar{\chi_{V_{\lambda}}}{P}_{\SSS_m}$ vanishes
for $\lambda\vdash m$ s.t. $\weight(\lambda)>\degw(P)$.
\item
\label{theo-degree-poly-biii}
Let $W_m$ be a representation of\/ $\SSS_m$. If $\rho_m(P)=\chi_{W_m}$, then:
$$\degw(P)\geq \weight(W_m)\,.$$
\end{enumerate}
\item
\label{theo-degree-poly-c}
Let $W_m=\bigoplus_{\lambda\vdash m}V_\lambda{}^{\ng_{\lambda}}$ be a representation of $\SSS_m$. The polynomial
$$\relax{\Xg_{W_m}:=\sumnl
_{\lambda\vdash m}\ng_{\lambda}
\Xg_\lambda\in\QQ[X_1,\ldots,X_{m-1}]}$$
where $\Xg_\lambda$ is the Frobenius 
polynomial for $\chi_{V_\lambda}$ (\ref{prop-Frobenius-char}), 
verifies
$$\rho_m(\Xg_{W_m})=\chi_{W_m}\text{\quad and\quad}
\degw(\Xg_{W_m})=\weight(W_m)\,,$$
and it is the only polynomial verifying simultaneously these two equalities, if and only if $\weight(W_m)\leq m/2$.
\end{enumerate}
\end{theo}
\demo
(\ref{theo-degree-poly-a}) In  \ref{Shur-ortho}, the inclusion `$\cont$' is \ref{prop-Frobenius-char}. For the converse, it suffices to restrict oneself to the case of a monomial $\rho_m(X_1^{\ng_1}\cdots X_d^{\ng_d})$ for $\sum_i i\,\ng_i\leq d$, or, thanks to \ref{prop-basic-CP}-(\ref{prop-basic-CP-b}), what amounts to the same thing, to show that
$$(\chi_{\IE_1^m})^{\ng_1}\cdots (\chi_{\IE_d^m})^{\ng_d}\in
\big\langle\chi_{V_\lambda}\mathrel{\big|}\big(|\lambda|=m\big)\ \&\ \big(\weight(\lambda)\leq d\,\big)\big\rangle{\vrule depth2pt width0pt}_{\QQ}\,.$$
Here, one recognizes, on the left, the character of $\IE_{\nu}^m$, for $\nu:=(1^{\ng_1},\ldots, d^{\ng_d})$. Now, since we already know after \ref{prop-basic-CP}-(\ref{prop-basic-CP-c}), that $\weight(\IE_{\nu}^m)\leq |\nu|$, we conclude that 
$\weight(\IE_{\nu}^m)\leq d$. 
The irreducible factors of $\IE_{\nu}^m$ are therefore of the form $V_\lambda$ with $\weight(\lambda)\leq d$, which settles the  inclusion `$\dans$'.

For the last claim about $\ker(\rho_{m|d})$, notice that the space $\QQ^{\leq d}[X_1,\ldots,X_d]$ admits as a basis the set of monomials 
$$\let\ag\ng
M(d):=\bigset X_1^{\ag_1}X_{2}^{\ag_2}\ldots X_{d}^{\ag_{d}}\mid\textstyle\sum_{i=1}^{d}i\,\ag_i\leq d\,/\,,$$
while the space $\QQ_{\class}^{\leq d}(\SSS_m)$ has a basis indexed by the following set of partitions\,:
$$\mathrigid3mu
L(d):=\bigset \ulambda:=(1^{\ng_1},2^{\ng_2},\ldots,\lambda_2^{\ng_{\lambda_2}})\mid
\big(\textstyle\sum_{i=1}^{\lambda_2}i\,\ng_i\leq d\big)\mathbin\&\big(\textstyle\sum_{i=1}^{\lambda_2}i\,\ng_i+\lambda_2\leq m\big)
/\,,$$
because $\set\chi_{V_{\lambda}}\mid\lambda\vdash m/$ is linearly independent after Schur orthogonality.

Since $\lambda_2\leq d$, the map $\xi:L(d)\to M(d)$, which associates $(1^{\ng_1},2^{\ng_2},\ldots,\lambda_2^{\ng_{\lambda_2}})$ with the monomial $X_1^{\ng_1}X_2^{\ng_2}\ldots X_{\lambda_2}^{\ng_{\lambda_2}}$ is well defined and injective. The condition for $\xi$ to be bijective is that $\lambda_2$ be able to take the value~$d$, in which case $\textstyle\sum_{i=1,\dots,\lambda_2}i\,\ng_i=d$, and that condition appears to be simply that $2d\leq m$.

\medskip
\noindent(\ref{theo-degree-poly-bi}) 
Again, thanks to \ref{prop-basic-CP}-(\ref{prop-basic-CP-b}), we need only prove that the number:
$$\bigscalar1{
\chiE_{1}^{m}(g)^{\ng_1}
\chiE_{2}^{m}(g)^{\ng_2}
\cdots
\chiE_{r}^{m}(g)^{\ng_r}}_{\SSS_m}\,,
\eqno(\diamond)$$
is constant for all 
$m\geq\sum\nolimits_{i=1}^{r} i\,\ng_i\,.$
But $(\diamond)$ 
is simply the dimension of the subspace 
of invariant tensors in the $\SSS_m$-module
$$\IE^{m}_{\nu}:=
(\IE_{1}^{m})^{\otimes\ng_1}\otimes
(\IE_{2}^{m})^{\otimes\ng_2}\otimes\cdots\otimes
(\IE_{r}^{m})^{\otimes\ng_r}\,,
$$
where $\nu:=(1^{\ng_1},2^{\ng_2},\ldots,r^{\ng_{r}})$.
Let $\nu=(\nu_1\geq\nu_2\geq\cdots\geq\nu_{\ell})$.
The canonical basis  $\B^{m}_{\nu}$ of $\IE_{\nu}^{m}$ is the set of tensors
$${\gamma_1}\otimes{\gamma_2}\otimes\cdots\otimes
{\gamma_\ell}\,,
$$
where $\gamma_i$ is a cycle of length $\nu_i$ of $\SSS_m$.
The basis $\B^{m}_\nu$ is clearly in bijection with the set $\T^{m}_{\nu}$ of Young \emph{tableaux\/} $\tau$ of shape $\nu$, with the $i$'th row  filled with elements of $\iii[1,m]$ 
so that they represent
a cycle of length $\nu_i$ of $\SSS_m$. The action of $\SSS_m$ on $\B^{m}_{\nu}$ induces in $\T^{m}_{\nu}$ the natural action of $\SSS_m$ on Young tableaux.

Because of these identifications, the dimension of $(\IE^{m}_{\nu}){}^{\SSS_m}$ is the cardinality of the orbit space $\T^{m}_\nu/\SSS_m$, and the stability we are seeking to prove is equivalent to the natural map
$$\T^{m}_\nu/\SSS_m\to\T^{m+1}_\nu/\SSS_{m+1}\,,\ [\tau\mod {\SSS_m}]\mapsto[\tau\mod {\SSS_{m+1}}]\,,$$
being a bijection. But the necessary and sufficient condition for this is precisely that the total number of boxes $|\nu|$ be smaller or equal to $m$, since, in that case, a single permutation $g\in \SSS_m$ will 
allow renumbering of all boxes simultaneously with numbers in the interval $\iii[1,|\nu|]$.
Hence $(\diamond)$ is constant for 
$m\geq |\nu|=\sumnl_{i}i\,\ng_i\,.$

\smallskip
\noindent(\ref{theo-degree-poly-bii}) After (\ref{theo-degree-poly-a}), $\rho_m(P)$ belongs to the linear span of the characters $\chi_{V_{\nu}}$ of $\SSS_m$, with $\weight(\nu)\leq\degw(P)$, and, thanks to Schur orthogonality, these characters are orthogonal to any $\chi_{V_{\lambda}}$ with $\weight(\lambda)> \degw(P)$, hence the claim.

\smallskip
\noindent(\ref{theo-degree-poly-biii}) Indeed, 
if we had $\degw(P)<\weight(W_m)$ then, for any irreducible factor $V_{\lambda}$ of $W_m$ of weight $\weight(V_{\lambda})=\weight(W_m)$, we would get, after (\ref{theo-degree-poly-bii})\,:
$$0=\scalar{\chi_{V_\lambda}}{P}_{\SSS_{m}}=\scalar{\chi_{V_\lambda}}{W_m}_{\SSS_{m}}\not=0\,,$$
which is a contradiction. Hence $\degw(P)\geq\weight(W_m)$.

\smallskip
\noindent(\ref{theo-degree-poly-c}) Immediate consequence of (\ref{theo-degree-poly-biii}) and the study of $\ker(\rho_{m|d})$ in (\ref{theo-degree-poly-a}).
\enddemo

\noendpoint\begin{comms}
\begin{itemize}
\item The proposition \ref{prop-basic-CP}-(\ref{prop-basic-CP-a}) showed that, for fixed $m\geq N$, the character of the tensor product $\IE^{m}_{\nu}:=
(\IE_{1}^{m})^{\otimes\ng_1}\otimes
(\IE_{2}^{m})^{\otimes\ng_2}\otimes\cdots\otimes
(\IE_{N}^{m})^{\otimes\ng_N}$, as an $\SSS_m$-module,
is the polynomial\,:
$\Eg_\nu:=\Eg_1^{\ng_1}\Eg_{2}^{\ng_{2}}\cdots \Eg_{N}^{\ng_{N}}$
of $\QQ[X_1,\ldots,X_N]\,,$
whose weight $\degw(\Eg_\nu)=\sum_i i\,\ng_i$ can be arbitrarily big. 
On the other hand, theorem \ref{theo-degree-poly}-(\ref{theo-degree-poly-c})
showed that the same character is given by the polynomial $\Xg_{\IE^{m}_{\nu}}$ of $\QQ[X_1,\ldots,X_{m-1}]$, of weight $\weight(\IE_{\nu}^{m})<m$.

This disagreement is due to the fact that, while $\Eg_\nu$ expresses a priori all the characters in the family 
$\set \chi_{\IE^{m}_{\nu}}/_{m\geq N}$ \emph{simultaneously}, the Frobenius polynomial $\Xg_{\IE^{m}_{\nu}}$ expresses a priori only one of them, viz. $\chi_{\IE^{m}_{\nu}}$. 

Heuristically speaking, the elements in $\ker(\rho_m)$ 
allow the reduction of the weight of the polynomial $\Eg_{\nu}$ in order to reach the required bound $m$.
This distinction will become meaningful in the next section.

\item It is also worth noting that while the explicit writing of $\Eg_{\nu}$ is very simple, thanks to formula \ref{prop-basic-CP}-(\ref{prop-basic-CP-a}), the writing of $\Xg_{\IE_{\nu}^{m}}$ can be 
quite complex insofar as it entails knowing the set of multiplicities
 $\ng_{\lambda}$ of the irreducible factors $V_{\lambda}$ in the decomposition $\IE_{\nu}^{m}=\bigoplus_{\lambda\vdash m} V_{\lambda}^{\ng_{\lambda}}$, and also in the explicit description of each $\Xg_{\lambda}$, for $\ng_{\lambda}\not=0$.
\end{itemize}
\end{comms}

\smallskip
\noindent The following is a useful corollary to theorem \ref{theo-degree-poly}-(\ref{theo-degree-poly-c})\killline

\begin{coro}\label{coro-degree-poly}Let $\W$ be an $\FB$-module which is  $\PC_{m\geq N}$ and has polynomial character $P_{\W}$.
\def\varlistskips{\itemsep2pt\mou\parsep0pt\mou\topsep0pt}\begin{enumerate}
\item
\label{coro-degree-poly-a}
$P_{\W}=\Xg_{W_m}$, for all $m\geq \max\set2\degw(P_{\W}),N/$.
\item
\label{coro-degree-poly-b}
 $\weight(W_m)=\degw(P_{\W})$, for all $m\geq  \max\set2\degw(P_{\W}),N/$.
\item
\label{coro-degree-poly-c}
 $\weight(\W_{\geq N})=\degw(P_{\W})$.
\end{enumerate}
\end{coro}
\demo (\ref{coro-degree-poly-a},\ref{coro-degree-poly-b}) Immediate after \ref{theo-degree-poly}-(\ref{theo-degree-poly-c}). \parskip1ex

\noindent
(\ref{coro-degree-poly-c}) For all $m\geq N$, we have $\rho_{m}(P_{\W})=\chi_{W_m}$ so that $\weight(W_m)\leq\degw(P_{\W})$, again
after \ref{theo-degree-poly}-(\ref{theo-degree-poly-c}). Therefore, 
$\sup_{m\geq N}\set \weight(W_m)/=\degw(P_{\W})$,
after (\ref{coro-degree-poly-b}).
\enddemo

\section{$\PC$ versus $\RS$}\label{PC-vs-RS}

\subsection{Equivalence of the stability conditions}
We show the equivalence, for an $\FB$-module $\W$, 
between the properties $\RS$ and $\PC$. In this, two more elements will deserve special attention. First, the ranks of validity of the properties, and, second, the weight of the polynomial~$P_{\W}$. The complete statement is the following.

\inhibe
\subsection{The $\FI$-modules $\EE_{1}^{\otimes\udp\ell}\dans\EE_{1}^{\otimes\ell}$}\label{El!}
For $m\geq1$, we endow the vector space $\IE_{1}^{m}=k^{m}$ with the action of $\SSS_m$ induced by its action on the canonical basis $\set e_i/_{i=1,\ldots,m}$ defined by $g(e_i):=e_{g(i)}$. 

For $\ell\geq1$ and any map $f:\iii[1,\ell]\to\iii[1,m]$, we denote by `$e_{f}$' the tensor\,:
$$e_f:=e_{f(1)}\otimes\cdots\otimes e_{f(\ell)}\in(\IE_{1})^{\otimes\ell}\,.$$
The set $\B^{\ell}:=\set e_{f}\mid f:\iii[1,\ell]\to\iii[1,m]/$
is a basis of $(\IE_{1}^{m})^{\otimes\ell}$. It will be endowed
 with the \emph{diagonal action} of $\SSS_m$, \idest 
$$g(e_{f}):=e_{g\circ f}\,,\mrlap{\quad\forall g\in\SSS_m\,.}$$
Now, let $\B^{\udp \ell}$ denote the subset of tensors $e_{f}$ with $f$ \emph{injective}. We have\,:
$$|\B^{\udp \ell}|=m^{\udp\ell}:=m(m-1)\cdots(m-\ell+1)\,,$$
and, in particular, $\B^{\udp\ell}=\emptyset$, if $\ell>m$.
The set $\B^{\udp \ell}$ is stabilized by the diagonal action of $\SSS_m$, and the subspace it spans, \idest
$$(\IE_{1}^{m})^{\otimes \udp\ell }:=\langle e_{f}\mid\text{$f$ is injective}\rangle_{k}\,,$$
is a sub-$\SSS_m$-module of $(\IE_{1}^{m})^{\otimes\ell}$.

\noendpoint\begin{prop}
\def\varlistskips{\itemsep4pt}\begin{enumerate}
\item The natural inclusion $(\IE_1^{m})^{\otimes\ell}\dans (\IE_1^{m+1})^{\otimes\ell}$ defining the $\FI$-module
$$\IE_{1}^{\otimes\ell}:=\set(\IE_1^{m})^{\otimes\ell}\to (\IE_1^{m+1})^{\otimes\ell}/_{m\geq0}\,,$$
restricts to an inclusion defining the sub-$\FI$-module
$(\IE_1^{m})^{\otimes\udp\ell}\dans (\IE_1^{m+1})^{\otimes\udp\ell}$
defining the sub-$\FI$-module
$$\IE_{1}^{\otimes\udp\ell}:=\set(\IE_1^{m})^{\otimes\udp\ell}\to (\IE_1^{m+1})^{\otimes\udp\ell}/_{m\geq0}\,.$$

\item The $\FI$-modules $\IE_{1}^{\otimes\ell}$ and $\IE_{1}^{\otimes\udp\ell}$ are $\PC_{m\geq0}$.

More precisely, the characters of their terms are respectively given by the polynomials\,:
$(X_1)^{\ell}$ and $(X_1)^{\udp\ell}$,
where $(X_1)^{\udp\ell}$ denotes the \emph{failing factorial\/} $X_1(X_1-1)\cdots(X_1-\ell+1)$.

\item The $\FI$-modules $\IE_{1}^{\otimes\ell}$ and $\IE_{1}^{\otimes\udp\ell}$ are $\RS_{m\geq 2\ell}$.
\end{enumerate}
\end{prop}

\subsection{}
\endinhibe

\begin{theo}\label{theo-PC=RS}Let $\W$ be an $\FB$-module.
\def\varlistskips{\advance\leftmargin1em}\begin{enumerate}
\item\label{theo-PC=RS-a} 
If $\W$ is $\RS$, then $\W$ is $\PC$ and $\rankPC(\W)\leq\rankRS(\W)$.

\item\label{theo-PC=RS-b}If $\W$ is $\PC$ with polynomial character $P_\W$, then $\W$ is $\RS$ and
$$
\rankRS(\W)\leq \max\set \rankPC(\W),2\degw(P_{\W})/\,.
$$
\end{enumerate}
\end{theo}
\demo (\ref{theo-PC=RS-a}) is corollary \ref{coro-Frobenius-char}.
(\ref{theo-PC=RS-b})
Set $d:=\degw(P_{\W})$. After \ref{coro-degree-poly}-(\ref{coro-degree-poly-a}),  we know that $P_{\W}=\Xg_{W_m}$ for all $m\geq2d$, so that, if we decompose $W_{2d}$ in its simple $\QQ[\SSS_{2d}]$-submodules\,:
$$W_{2d}=\bigoplusnl_{|\lambda|\leq d}(V_{\lambda[2d]})^{\ng_{\lambda}}\,,$$
and if we denote 
$\lambda':=\parti(\lambda_1,\lambda_1,\lambda_2,\ldots,\lambda_\ell)$
for $\lambda:=\parti(\lambda_1,\lambda_2,\ldots,\lambda_\ell)$, then
the $\FI$-module
$$\W':=\bigoplusnl_{|\lambda|\leq d}(\V_{\lambda'})^{\ng_{\lambda}}\,,
$$
which is clearly $\RS_{m\geq2d}$, will be such that 
$P_{\W'}=P_{\W}$, by construction.
As a consequence, there exists an isomorphism of $\FB$-modules:
$$\W_{\geq\max\set 2d,N/}\sim\W'_{\geq \max\set 2d,N/}\,,\eqno(\diamond)$$
and $\W$ is $\RS$ for $m\geq \max\set 2d,N/$, as stated.
\enddemo

\begin{rema}\label{last-rema}The last proof shows, in $(\diamond)$, that an $\FB$-module that is $\PC$ is asymptotically isomorphic, as $\FB$-module, to an $\RS$ $\FI$-module.
\end{rema}

\section{Addendum on the weight of a tensor product}
\vskip-0.5ex\vskip0pt
\subsection{Reinterpretation of the weight of a representation}The theorem
\ref{theo-degree-poly}-(\ref{theo-degree-poly-c}) gives an alternative definition of the weight of a representation $W_m$ of $\SSS_m$ \emph{as the lowest possible weight of the polynomials expressing its character $\chi_{W_m}$}. \killline

\subsection{On the weight of the tensor product of representations}

\noendpoint\begin{prop}\label{prop-weight-tensor-product}
\def\varlistskips{\itemsep2pt\mou\parsep0pt\mou\topsep2pt}
\begin{enumerate}
\item
\label{prop-weight-tensor-product-a}
If $W_{1,m}$ and $W_{2,m}$ are representations of $\SSS_m$, then
$$\weight(W_{1,m}\otimes W_{2,m})\leq\weight(W_{1,m})+\weight(W_{2,m})\,.$$
\item
\label{prop-weight-tensor-product-b}
If $\W_1$ and $\W_2$ are $\FB$-modules which are $\PC_{m\geq N}$, then 
$$\weight(W_{1,m}\otimes W_{2,m})=\weight(W_{1,m})+\weight(W_{2,m})\,,$$
for all $m\geq 2(\weight(W_{1,m})+\weight(W_{2,m}))$. In particular,
$$\weight((\W_{1}\otimes\W_{2})_{\geq N})=
\weight((\W_{1})_{\geq N})+
\weight((\W_{2})_{\geq N})
$$
\item
\label{prop-weight-tensor-product-c}
\hfil $\weight(V_{\lambda[m]}\otimes V_{\mu[m]})=
|\lambda|+|\mu|$,\quad $\forall m\geq 2(|\lambda|+|\mu|)$.
\end{enumerate}
\end{prop}
\demo (\ref{prop-weight-tensor-product-a})\parskip3pt\mou
 As $\chi_{W_1\otimes W_2}=\chi_{W_1}\chi_{W_2}$, we get
%
$$\mathrigid1.2mu
\weight(W_1\otimes W_2)\mathop\leq_{(1)}\degw(\Xg_{W_1}\Xg_{W_2})
=\degw(\Xg_{W_1})+\degw(\Xg_{W_2})\mathop=_{(2)}\weight(W_1)+\weight(W_2)\,,
$$
inequality `$\,\leq_{(1)}$' after \ref{theo-degree-poly}-(\ref{theo-degree-poly-biii}), and equality `$\,=_{(2)}$' after \ref{theo-degree-poly}-(\ref{theo-degree-poly-c}).

\smallskip
\noindent (\ref{prop-weight-tensor-product-b}) Let $d_i:=\degw(P_{\W_i})$. We know after corollary \ref{coro-degree-poly}-(\ref{coro-degree-poly-a}), that 
$P_{\W_i}=\Xg_{W_{i,m}}\,,$
for $m\geq 2\max\set d_1,d_2/$. Hence,
$\Xg_{W_{1,m}}\Xg_{W_{2,m}}=\Xg_{W_{1,m}\otimes W_{2,m}}\,,$
for $m\geq 2(d_1+d_2)$, after the same corollary. Therefore,
$$\weight(W_{1,m}\otimes W_{2,m})=\weight (W_{1,m})+\weight (W_{2,m})\,,
$$
for $m\geq 2(d_1+d_2)$, again after \ref{theo-degree-poly}-(\ref{theo-degree-poly-c}).

\smallskip
\noindent (\ref{prop-weight-tensor-product-c}) Particular case of 
(\ref{prop-weight-tensor-product-b}). 
\enddemo



\end{document}